\documentclass{amsart}
\usepackage{amssymb, latexsym}
\usepackage{enumerate}

\newcommand\card{\operatorname{card}}

\theoremstyle{plain}
\newtheorem{thm}{Theorem}
\newtheorem{conjecture}{Conjecture}
\newtheorem{condition}[thm]{Condition}
\newtheorem{rem}[thm]{Remark}
\newtheorem{prop}{Proposition}
\newtheorem{lem}[prop]{Lemma}

\numberwithin{equation}{section} \numberwithin{thm}{section}

\begin{document}

\title[Barely $\dot{H}^{s_{p}}$ Supercritical Wave Equation ] { One Remark On Barely $\dot{H}^{s_{p}}$ Supercritical
Wave Equations}
\author{Tristan Roy}
\address{University of California, Los Angeles}
\email{triroy@math.ucla.edu}

\begin{abstract}
We prove that a good $\dot{H}^{s_{p}}$ critical theory for the $3D$
wave equation $\partial_{tt} u - \triangle u = -|u|^{p-1} u$ can be
extended to prove global well-posedness of smooth solutions of at
least one $3D$ barely $\dot{H}^{s_{p}}$ supercritical wave equation
$\partial_{tt} u - \triangle u =- |u|^{p-1} u g(|u|)$, with $g$
growing slowly to infinity, provided that a Kenig-Merle type
condition is satisfied. This result extends those
\cite{taolog,triroy} obtained for the particular case $s_{p}=1$.

 \end{abstract}

\maketitle

\section{Introduction}

We shall consider the following wave equation

\begin{equation}
\left\{
\begin{array}{ll}
\partial_{tt} u - \triangle u & = -|u|^{p-1} u g \left(|u| \right)
\\ & \\
u(0) & :=u_{0} \in \tilde{H}^{2} \\ & \\
\partial_{t} u(0) & := u_{1} \in \tilde{H}^{1}
\end{array}
\right. \label{Eqn:BarelyWaveEq}
\end{equation}
where $u: \mathbb{R} \times \mathbb{R}^{3} \rightarrow \mathbb{C}$
is a complex-valued scalar field, $p>3$, $\tilde{H}^{2}:=
\dot{H}^{2}(\mathbb{R}^{3}) \cap \dot{H}^{s_{p}}(\mathbb{R}^{3})$,
$\tilde{H}^{1}:=\dot{H}^{1}(\mathbb{R}^{3}) \cap
\dot{H}^{s_{p}-1}(\mathbb{R}^{3})$ and $g$ is a smooth, real-valued
positive function defined on the set of nonnegative numbers and
satisfying

\begin{equation}
\begin{array}{lll}
0 \leq & g^{'}(x) & \lesssim  \frac{1}{x}
\end{array}
\label{Eqn:Conditiong2}
\end{equation}
Condition (\ref{Eqn:Conditiong2}) basically says that $g$ grows more
slowly than any positive power of $u$.

We shall see that this equation (\ref{Eqn:BarelyWaveEq}) has many
connections with the defocusing power-type defocusing wave equation

\begin{equation}
\left\{
\begin{array}{ll}
\partial_{tt} u - \triangle u & = - |u|^{p-1} u \\ & \\
u(0) & := u_{0} \in \dot{H}^{s_{p}} (\mathbb{R}^{3}) \\ & \\
\partial_{t} u(0) & := u_{1} \in \dot{H}^{s_{p}-1} (\mathbb{R}^{3})
\end{array}
\right. \label{Eqn:PDef}
\end{equation}
Here $s_{p}:= \frac{3}{2} - \frac{2}{p-1}$. It is known that if $u$
satisfies (\ref{Eqn:PDef}) then $u_{\lambda}$ defined by

\begin{equation}
\begin{array}{ll}
u_{\lambda}(t,x) & := \frac{1}{\lambda^{\frac{2}{p-1}}} u \left(
\frac{t}{\lambda}, \, \frac{x}{\lambda} \right)
\end{array}
\label{Eqn:ScalingEq}
\end{equation}
satisfies the same equation, but with data $ u_{\lambda}(0,x)=
\frac{1}{\lambda^{\frac{2}{p-1}}} u_{0} \left( \frac{x}{\lambda}
\right) $ and $ \partial_{t}u_{\lambda}(0,x)=
\frac{1}{\lambda^{\frac{2}{p-1}+1}} u_{1} \left( \frac{x}{\lambda}
\right) $. Notice that (\ref{Eqn:PDef}) is
$\dot{H}^{s_{p}}(\mathbb{R}^{3})$ critical, which means that the
$\dot{H}^{s_{p}}(\mathbb{R}^{3}) \times
\dot{H}^{s_{p}-1}(\mathbb{R}^{3})$ norm of $\left( u(0),
\partial_{t} u(0) \right)$ is invariant under the scaling defined
above.

We recall the local existence theory: it is known (see for example
\cite{kenigmerle1,sogge}) that there exists a  positive constant
$\delta:=\delta \left( \| ( u_{0}, u_{1}) \|_{\dot{H}^{s_{p}}
(\mathbb{R}^{3}) \times \dot{H}^{s_{p}-1}(\mathbb{R}^{3})} \right)
> 0$, a time of local existence $T_{l}>0$ such that if

\begin{equation}
\begin{array}{ll}
\| \cos{(tD)} u_{0} + \frac{\sin{(tD)}}{D} \|_{ L_{t}^{2(p-1)}
L_{x}^{2(p-1)} ( [0,T_{l}] \times \mathbb{R}^{3} ) } & \leq  \delta
\end{array}
\end{equation}
then there exists a unique solution $ \left( u,
\partial_{t} u \right) \in \mathcal{C} \left( [0,T_{l}], \, \dot{H}^{s_{p}}(\mathbb{R}^{3}) \right)
\cap L_{t}^{2(p-1)} L_{x}^{2(p-1)}( [0,T_{l}] \times \mathbb{R}^{3})
\cap D^{\frac{1}{2} - s_{p}} L_{t}^{4} L_{x}^{4} ( [0,T_{l}] \times
\mathbb{R}^{3}) \times \mathcal{C} \left( [0,T_{l}], \,
\dot{H}^{s_{p}-1} (\mathbb{R}^{3}) \right) $ of (\ref{Eqn:PDef})
\footnote{Notice that the $L_{t}^{2(p-1)} L_{x}^{2(p-1)}(\mathbb{R}
\times \mathbb{R}^{3})$ norm of $u$ is invariant under the scaling
(\ref{Eqn:ScalingEq}). The choice of the space $ L_{t}^{2(p-1)}
L_{x}^{2(p-1)}$ in which we place the solution $u$ is not unique.
There exists un infinite number of spaces of the form $L_{t}^{q}
L_{x}^{r}$ scale invariant in which we can establish a local
well-posedness theory} in the sense of the integral equation, i.e
$u$ satisfies the Duhamel formula

\begin{equation}
\begin{array}{ll}
u(t) & : = \cos{(tD)} u_{0} + \frac{\sin{(tD)}}{D} u_{1} -
\int_{0}^{t} \frac{\sin{(t-t^{'}) D}}{D} \left( |u|^{p-1} u \right)
(t^{'}) \, dt^{'}
\end{array}
\label{Eqn:DuhamelCritical}
\end{equation}
It follows that we can define a maximal time interval of existence
$I_{max}=(-T_{-},T_{+})$. Moreover $\| u \|_{L_{t}^{2(p-1)}
L_{x}^{2(p-1)}(J)} < \infty$, $\| D^{s_{p}-\frac{1}{2}} u
\|_{L_{t}^{4} L_{x}^{4}(J) } < \infty $ and \\ $ \| (u,\partial_{t}
u) \|_{L_{t}^{\infty} \dot{H}^{s_{p}} \times L_{t}^{\infty}
\dot{H}^{s_{p}-1} (J)} < \infty $ for all subinterval $J \subset
I_{max}$. See \cite{kenigmerleenergy} or \cite{taobook} for more
explanations.

Now we turn to the global well-posedness theory of (\ref{Eqn:PDef}).
In view of the local well-posedness theory, one can prove (see
\cite{kenigmerle1} and references), after some effort, that it is
enough to find a finite upper bound of $\| u \|_{L_{t}^{2(p-1)}
L_{x}^{2(p-1)}(I \times \mathbb{R}^{3})}  $ on arbitrary long time
intervals $I$, and, if this is the case, then the solution scatters
to a solution of the linear wave equation. No blow-up has been
observed for (\ref{Eqn:PDef}). Therefore it is believed that the
following scattering conjecture is true

\begin{conjecture}{"\textbf{Scattering Conjecture}"}
Assume that $u$ is the solution of (\ref{Eqn:PDef}) with data
$(u_{0},u_{1}) \in \dot{H}^{s_{p}}(\mathbb{R}^{3}) \times
\dot{H}^{s_{p}-1}(\mathbb{R}^{3})$. Then $u$ exists for all time $t$
and there exists $C_{1}:=C_{1} \left( \| ( u_{0} , u_{1})
\|_{\dot{H}^{s_{p}}(\mathbb{R}^{3}) \times
\dot{H}^{s_{p}-1}(\mathbb{R}^{3}) } \right) $ such that
\begin{equation}
\begin{array}{ll}
\| u \|_{L_{t}^{2(p-1)} L_{x}^{2(p-1)}(\mathbb{R} \times
\mathbb{R}^{3})} & \leq C_{1}
\end{array}
\end{equation}
\label{Conj:1}
\end{conjecture}

The case $s_{p}=1$ (or, equivalently, $p=5$) is particular. Indeed
the solution $(u,\partial_{t}u) \in \mathcal{C} \left( [0,T_{l}], \,
\dot{H}^{1}(\mathbb{R}^{3}) \right) \times \mathcal{C} \left(
[0,T_{l}], \, L^{2} (\mathbb{R}^{3} ) \right)$ satisfies the
conservation of the energy $E(t)$ defined by

\begin{equation}
\begin{array}{ll}
E(t) & : = \frac{1}{2} \int_{\mathbb{R}^{3}}  | \partial_{t} u (t,x)
|^{2} \, dx + \frac{1}{2} \int_{\mathbb{R}^{3}} |\nabla u (t,x)|^{2}
\, dx + \frac{1}{6} \int_{\mathbb{R}^{3}} |u|^{6}(t,x) \, dx
\end{array}
\label{Eqn:DefNrj}
\end{equation}
In other words, $E(t)=E(0)$. This is why this equation is often called energy-critical: the exponent $s_{p}=1$ precisely corresponds to the
minimal regularity required for (\ref{Eqn:DefNrj}) to be satisfied. The global well-posedness of (\ref{Eqn:DefNrj}) in the energy class and in
higher regularity spaces is now understood. Rauch \cite{rauch} proved the global existence of smooth solutions of this equation with small data.
Struwe \cite{St} showed that the result still holds for large data but with the additional assumption of spherical symmetry of the data. The
general case (large data, no symmetry assumption) was finally settled by Grillakis \cite{grill,grill2}. Later Shatah and Struwe
\cite{shatstruwe} reproved this result. Kapitanski \cite{kap} and, independently, Shatah and Struwe \cite{shatstruwe2}, proved global existence
of solutions in the energy class. Bahouri and Gerard \cite{bahger} reproved this result by using a compactness method. In particular, they
showed that the $L_{t}^{2(5-1)} L_{x}^{2(5-1)} (\mathbb{R} \times \mathbb{R}^{3})$ norm of the solution is bounded by an unspecified finite
quantity. Lately Tao \cite{tao} found an exponential tower type bound of this norm. Notice that, in all these proofs of global existence of
solutions of the energy-critical wave equation, the conservation of energy, which leads, in particular, to the control of the $\dot{H}^{1}
\times L^{2}$ norm of the solution $\left( \partial_{t} u(t), u(t) \right)$, is a key point.

If $s_{p}<1$, or equivalently, $p<5$, then we are in the
energy-subcritical equation. The scattering conjecture is an open
problem. Nevertheless, some partial results are known if we consider
the same problem (\ref{Eqn:PDef}), but with data $(u_{0},u_{1}) \in
H^{s} \times H^{s-1}$, $ s_{p} < s $. More precisely, it is proved
in \cite{kpv2,triroysub} that there exists $s_{0}:=s_{0}(p)$ such
that $s_{p} < s_{0} < 1$ and such that (\ref{Eqn:PDef}) is globally
well-posed in $H^{s} \times H^{s-1}$, for $s>s_{0}$. These results
rely upon two well-known strategies: the Fourier truncation method
\cite{kpv2}, designed by Bourgain \cite{bourg,bourg2}, and the
$I$-method \cite{triroysub}, designed by
Colliander-Keel-Staffilani-Takaoka-Tao \cite{almckstt}.

If $s_{p}>1$, or, equivalently, $p>5$, then we are in the
energy-supercritical regime. The global behaviour of the solution
is, in this regime, very poorly understood. Indeed, if we follow the
theory of the energy-critical wave equation, then the first step
would be to prove that the $\dot{H}^{s_{p}} \times
\dot{H}^{s_{p}-1}$ norm of the solution is bounded for all time by a
finite quantity that would only depend on the $\dot{H}^{s_{p}}
\times \dot{H}^{s_{p}-1}$ norm of the initial data. Unfortunately,
the control of this norm is known to be a very challenging problem,
since there are not known conservation laws in high regularity
Sobolev spaces. Lately Kenig and Merle \cite{kenigmerle1} proved, at
least for radial data, that this step would be the last, by using
their concentration compactness/rigidity theorem method
\cite{kenigmerleenergy}. More precisely, they showed that if $
\sup_{t \in I_{max}} \| ( u(t) ,\partial_{t} u(t))
\|_{\dot{H}^{s_{p}} (\mathbb{R}^{3}) \times
\dot{H}^{s_{p}-1}(\mathbb{R}^{3})} < \infty $, then Conjecture
\ref{Conj:1} is true. See also the recent work of Killip-Visan
\cite{kilvis}, which establishes results in this direction for the
Schr\"odinger equation.

As mentioned before, the energy supercritical regime is almost
\textit{terra incognita}. Nevertheless, Tao \cite{taolog} observed
that the technology used to prove global well-posedness of smooth
solutions of (\ref{Eqn:PDef}) can be extended, after some effort, to
some equations of the type (\ref{Eqn:BarelyWaveEq}), with  $p=5$ and
radial data. More precisely, he proved global regularity of
(\ref{Eqn:BarelyWaveEq}) with $ g(x)  :=\log{(2+ x^{2})} $. This
phenomenon, in fact, does not depend on the symmetry of the data: it
was proved in \cite{triroy} that there exists a unique global smooth
solution of (\ref{Eqn:BarelyWaveEq}) with $ g(x) := \log^{c}\log{(10
+ x^{2})} $ and $ 0 < c < \frac{8}{225}$.

Equations of the type (\ref{Eqn:BarelyWaveEq}) are called barely
$\dot{H^{s_{p}}}$-supercritical wave equations. Indeed, Condition
(\ref{Eqn:Conditiong2}) basically says that for every $\epsilon >
0$, there exist two constants $c_{1}:=c_{1}(p)$ and
$c_{2}:=c_{2}(p,\epsilon)$ such that for $|u|$ large

\begin{equation}
c_{1}(p) \leq  g(|u|) \leq c_{2}(p,\epsilon) |u|^{\epsilon}
\label{Eqn:Conditiong5}
\end{equation}
Since the critical exponent of the equation $
\partial_{tt} u - \triangle u = - |u|^{p-1 +  \epsilon} u $ is
$s_{p+ \epsilon}= s_{p} + O(\epsilon) $, the nonlinearity of
(\ref{Eqn:BarelyWaveEq}) is barely $\dot{H}^{s_{p}}$ supercritical.

The goal of this paper is to check that this phenomenon, observed
for $s_{p}=1$, still holds for other values of $s_{p}$. The standard
local well-posedness theory shows us that it is enough to control
the pointwise-in-time $\tilde{H}^{2} \times \tilde{H}^{1}$ norm of
the solution. In this paper, we will use an alternative local
well-posedness theory. We shall prove the following proposition

\begin{prop}{\textbf{"Local Existence for barely $\dot{H}^{s_{p}}$ supercritical wave equation"}}
Assume that $g$ satisfies  (\ref{Eqn:Conditiong2}) and

\begin{equation}
\begin{array}{ll}
g^{''}(x) & =O \left( \frac{1}{ x^{2}} \right)
\end{array}
\label{Eqn:Conditiong3}
\end{equation}
Let $M$ be such that $ \| (u_{0}, \ u_{1}) \|_{ \tilde{H}^{2} \times
\tilde{H}^{1}} \leq M$. Then there exists $\delta:= \delta (M)
> 0$ small  such that if $T_{l}$ satisfies

\begin{equation}
\begin{array}{ll}
 \| \cos{(tD)} u_{0} + \frac{\sin{tD}}{D} u_{1} \|_{L_{t}^{2(p-1)} L_{x}^{2(p-1)}([0,T_{l}] \times
 \mathbb{R}^{3})} \leq \delta
\end{array}
\label{Eqn:SmallCond}
\end{equation}
then there exists a unique $(u,\partial_{t} u ) \in \mathcal{C}
\left( [0, \, T_{l}], \, \tilde{H}^{2} \right) \cap L_{t}^{2(p-1)}
L_{x}^{2(p-1)} \left(  [0, \, T_{l}] \right) \cap D^{\frac{1}{2}
-s_{p}} L_{t}^{4} L_{x}^{4} \left(  [0, \, T_{l}] \right) \cap
D^{\frac{1}{2} -2} L_{t}^{4} L_{x}^{4} \left( [0, \, T_{l}] \right)
\times \mathcal{C} \left( [0, \, T_{l}], \, \tilde{H}^{1} \right)$
of (\ref{Eqn:BarelyWaveEq}) in the sense of the integral equation,
i.e $u$ satisfies the following Duhamel formula:

\begin{equation}
\begin{array}{ll}
u(t) & : = \cos{(tD)} u_{0} + \frac{\sin{tD}}{D} u_{1} -
\int_{0}^{t} \frac{\sin(t-t^{'})D}{D} \left( |u(t^{'})|^{p-1}
u(t^{'}) g(|u(t^{'})|) \right) \, dt^{'}
 \end{array}
\label{Eqn:Duhamelg}
\end{equation}
\label{prop:localexistencebar}
\end{prop}
Notice that there are many similarities between Proposition
\ref{prop:localexistencebar} and the local well-posedness theory for
(\ref{Eqn:PDef}).

This allows to define a maximum time interval of existence
$I_{max,g}=[-T_{-,g}, \, T_{+,g}]$ such that for all $J \subset
I_{max,g}$, we have $ \| u \|_{ L_{t}^{2(p-1)} L_{x}^{2(p-1)}(J) } <
\infty $,  $ \| D^{s_{p} -\frac{1}{2}} u \|_{W(J)} < \infty  $,  $
\| D^{2 -\frac{1}{2}} u \|_{W(J)} < \infty $ and $ \|(u,\partial_{t}
u) \|_{L_{t}^{\infty} \tilde{H}^{2}(J) \times L_{t}^{\infty}
\tilde{H}^{1}(J)} < \infty $. Again, see \cite{kenigmerleenergy} or
\cite{taobook} for more explanations.

Now we set up the problem. In view of the comments above for
$s_{p}=1$, we need to make two assumptions. First we will work with
a $`` \mathrm{good} ''$ $\dot{H}^{s_{p}}(\mathbb{R}^{3})$ theory:
therefore we will assume that Conjecture \ref{Conj:1} is true. Then,
we also would like to work with $\dot{H}^{s_{p}}(\mathbb{R}^{3})
\times \dot{H}^{s_{p}-1}(\mathbb{R}^{3})$ bounded solutions
$(u(t),\partial_{t} u(t))$; more precisely we will assume that the
following Kenig-Merle type condition holds

\begin{condition}{"\textbf{Kenig-Merle type condition}"}
Let $g$ be a function that satisfies (\ref{Eqn:Conditiong2}) and
that is constant for $x$ large. Then there exists $C_{2}:=C_{2}
\left( \| (u_{0}, \, u_{1}) \| _{\tilde{H}^{2} \times
\tilde{H}^{1}}, \, g \right)$ such that

\begin{equation}
\begin{array}{ll}
\sup_{t \in I_{max,g}} \| \left( u(t),\partial_{t} u(t) \right)
\|_{\dot{H}^{s_{p}}(\mathbb{R}^{3}) \times
\dot{H}^{s_{p}-1}(\mathbb{R}^{3})} & \leq C_{2}
\end{array}
\end{equation}
\label{Conj:2}
\end{condition}

\begin{rem}
In the particular case $s_{p}=1$, it is not difficult to see that
Condition \ref{Conj:2} is satisfied. Indeed $u$ satisfies the energy
conservation law

\begin{equation}
\begin{array}{ll}
E_{b}(t) & := \frac{1}{2} \int_{\mathbb{R}^{3}}  \left(
\partial_{t} u (t,x) \right)^{2} \, dx + \frac{1}{2}
\int_{\mathbb{R}^{3}} |\nabla u (t,x)|^{2} \, dx +
\int_{\mathbb{R}^{3}} F(u(t,x), \bar{u}(t,x)) \, dx
\end{array}
\end{equation}
with

\begin{equation}
\begin{array}{ll}
F(z,\bar{z}) & = |z|^{5+1} \int_{0}^{1} t^{5} \, \Re \left( g( t|z|
) \right) \, dt \\
& = |z|^{5+1} \int_{0}^{1} t^{5} \, g( t|z| )  \, dt
\end{array}
\end{equation}
Since $g$ is bounded then  $\left| F(z, \bar{z}) \right| \lesssim
|z|^{6}$. By using the Sobolev embeddings $ \| u_{0} \|_{L_{x}^{6}}
\lesssim \| u_{0} \|_{\tilde{H}^{2}}$ and $ \| u(t) \|_{L_{x}^{6}}
\lesssim \| u(t) \|_{\tilde{H}^{2}}$, we easily conclude that
Condition \ref{Conj:2} holds. The energy conservation law was
constantly used in \cite{taolog,triroy}.
\end{rem}
The main result of this paper is

\begin{thm}
There exists a function $\tilde{g}$ satisfying
(\ref{Eqn:Conditiong2}) and such that

\begin{equation}
\begin{array}{ll}
\lim_{x \rightarrow \infty} \tilde{g}(x) & = \infty
\end{array}
\end{equation}
such that the solution of (\ref{Eqn:BarelyWaveEq}) (with $g:=
\tilde{g}$) exists for all time, provided that the scattering
conjecture and the Kenig-Merle type condition are satisfied.
Moreover there exists a function $f$ depending on $T$ and $ \| (
u_{0}, u_{1}) \|_{\tilde{H}^{2} \times \tilde{H}^{1}}$ such that

\begin{equation}
\begin{array}{ll}
\| u \|_{L_{t}^{\infty} \tilde{H}^{2}([-T,T])}   + \| \partial_{t} u
\|_{L_{t}^{\infty} \tilde{H}^{1}([-T,T])}  & \leq f \left( T, \|(
u_{0},u_{1}) \|_{\tilde{H}^{2} \times \tilde{H}^{1}} \right)
\end{array}
\label{Eqn:BounduLinftyf}
\end{equation}
\label{Thm:GlobWave}
\end{thm}
Theorem \ref{Thm:GlobWave} shows that a ``good'' $
\dot{H}^{s_{p}}(\mathbb{R}^{3}) $ theory for (\ref{Eqn:PDef}) can be
extended, at least, to one barely $\dot{H}^{s_{p}}(\mathbb{R}^{3})$
supercritical equation, with $\tilde{g}$ going to infinity.

\begin{rem}
\begin{itemize}

\item $\tilde{g}$ is ``universal'': it does not depend on an upper bound of
the initial data

\item $\tilde{g}$ is unbounded: it goes to infinity as $x$ goes to infinity
\end{itemize}

\end{rem}

\begin{rem}
In fact, Theorem \ref{Thm:GlobWave} also holds for a weaker version
of Condition \ref{Conj:2}: there exists a function $C_{2}$ such that
for all subinterval $I \subset I_{max,g}$

\begin{equation}
\begin{array}{ll}
\sup_{t \in I} \left\| (u(t),\partial_{t}u(t) )
\right\|_{\dot{H}^{s_{p}} (\mathbb{R}^{3}) \times \dot{H}^{s_{p}-1}
(\mathbb{R}^{3}) } & \leq C_{2}
\end{array}
\end{equation}
with $C_{2}:=C_{2} \left( \| (u_{0},u_{1}) \|_{\tilde{H}^{2} \times
\tilde{H}^{1}}, \, g, \, |I| \right)$. See the proof of Theorem
\ref{Thm:GlobWave} and, in particular, (\ref{Eqn:ControlHsp}),
(\ref{Eqn:ControlHsp1}) and (\ref{Eqn:Controlnormi}).

\end{rem}

We recall some basic properties and estimates. If $t_{0} \in
[t_{1},t_{2}]$, if $F \in L_{t}^{\tilde{q}} L_{x}^{ \tilde{r}}
([t_{1}, \, t_{2}])$ and if  $( u ,\partial_{t} u) \in C \left(
[t_{1},t_{2}], \, \dot{H}^{m} (\mathbb{R}^{3}) \right) \times C
\left( [t_{1},t_{2}], \, \dot{H}^{m-1} (\mathbb{R}^{3}) \right)$
satisfy
\begin{equation}
\begin{array}{ll}
u(t) & : \cos{(tD)} u_{0} + \frac{\sin{tD}}{D} u_{1} -
\int_{t_{0}}^{t} \frac{\sin{(t-t^{'}) D}}{ D} F (t^{'}) \, d t^{'}
\end{array}
\label{Eqn:DuhForm}
\end{equation}
with data $(u(t_{0}),\partial_{t}u(t_{0})) \in
\dot{H}^{m}(\mathbb{R}^{3}) \times \dot{H}^{m-1}(\mathbb{R}^{3})$
then we have the Strichartz estimates \cite{ginebvelo,linsog}

\begin{equation}
\begin{array}{ll}
\| u \|_{L_{t}^{q} L_{x}^{r} ( [t_{1}, \, t_{2}] )} + \| u
\|_{L_{t}^{\infty} \dot{H}^{m}(\mathbb{R}^{3})  \left( [t_{1}, \,
t_{2}] \right)} + \| \partial_{t} u \|_{L_{t}^{\infty}
\dot{H}^{m-1}(\mathbb{R}^{3}) ( [t_{1}, \, t_{2}] )} & \lesssim
\left\| ( u(t_{0}), \,
\partial_{t} u(t_{0})) \right\|_{\dot{H}^{m}(\mathbb{R}^{3}) \times
\dot{H}^{m-1}(\mathbb{R}^{3}) } \\
& + \| F \|_{L_{t}^{\tilde{q}} L_{x}^{\tilde{r}} ( [t_{1}, \,
t_{2}])}
\end{array}
\label{Eqn:Strichartz}
\end{equation}
Here

\begin{itemize}

\item $(q,r)$ is $m$- wave admissible, i.e

\begin{equation}
\left\{
\begin{array}{l}
(q,r) \in (2, \infty] \times [2,\infty] \\
\frac{1}{q} + \frac{3}{r} = \frac{3}{2} - m \\
\end{array}
\right.
\end{equation}

\item

\begin{equation}
\begin{array}{l}
\frac{1}{q} + \frac{3}{r} = \frac{1}{\tilde{q}} + \frac{3}{
\tilde{r}} -2
\end{array}
\end{equation}
\end{itemize}

We set some notation that appear throughout the paper.

We write $A \lesssim B$ if there exists a universal nonnegative
constant $C^{'}>0$ such that $A \leq C^{'} B$. $A=O(B)$ means $A
\lesssim B$. More generally we write $A \lesssim_{a_{1},....,a_{n}}
B$ if there exists a nonnegative constant $C^{'}=C(a_{1},...,a_{n})$
such that $A \leq C^{'} B$. We say that $C^{''}$ is the constant
determined by $\lesssim$ in $A \lesssim_{a_{1},...,a_{n}} B $ if
$C^{''}$ is the smallest constant among the $C^{'}$ s such that $A
\leq C^{'}B$. We write $A <<_{a_{1},..,a_{n}} B$ if there exists a
universal nonnegative small constant $c=c(a_{1},...,a_{n})$ such
that $A \leq c B$.
Following \cite{kenigmerle1}, we define, on an
interval $I$

\begin{equation}
\begin{array}{ll}
\| u \|_{S(I)} & := \| u \|_{L_{t}^{2(p-1)} L_{x}^{2(p-1)} (I)} \\
\| u \|_{W(I)} & := \| u \|_{L_{t}^{4} L_{x}^{4}(I)} \\
\| u \|_{\tilde{W}(I)} & := \| u \|_{L_{t}^{\frac{4}{3}}
L_{x}^{\frac{4}{3}}(I)}
\end{array}
\end{equation}
We also define the following quantity

\begin{equation}
\begin{array}{ll}
Q(I,u) & := \| D^{s_{p} -\frac{1}{2}} u \|_{W(I)} + \|
D^{2-\frac{1}{2}} u \|_{W(I)}  + \| u \|_{L_{t}^{\infty}
\tilde{H}^{2}(I)} + \| \partial_{t} u \|_{L_{t}^{\infty}
\tilde{H}^{1}(I)}
\end{array}
\end{equation}
Let $X$ be a Banach space and $r \geq 0$. Then

\begin{equation}
\begin{array}{ll}
\mathbf{B}(X, r) & := \left\{ f \in X, \, \| f \|_{X} \leq r
\right\}
\end{array}
\end{equation}
We recall also the well-known Sobolev embeddings. We have

\begin{equation}
\begin{array}{ll}
\| h \|_{L^{\infty}(\mathbb{R}^{3})} & \lesssim \| h
\|_{\tilde{H}^{2}}
\end{array}
\label{Eqn:SobEmbed1}
\end{equation}
and

\begin{equation}
\begin{array}{ll}
\| h \|_{S(I)} & \lesssim \| D^{s_{p} -\frac{1}{2}} h
\|_{L_{t}^{2(p-1)} L_{x}^{\frac{6(p-1)}{2p-3}}(I)}
\end{array}
\label{Eqn:SobEmbed2}
\end{equation}
We shall combine (\ref{Eqn:SobEmbed2}) with the Strichartz
estimates, since  $ \left( 2(p-1), \, \frac{6(p-1)}{2p-3} \right)$
is $\frac{1}{2}$- wave admissible.

We also recall some Leibnitz rules \cite{christwein,kpv}. We have

\begin{equation}
\begin{array}{ll}
\| D^{\alpha} F(u) \|_{L_{t}^{q} L_{x}^{r} (I)} & \lesssim \|
F^{'}(u) \|_{L_{t}^{q_{1}} L_{x}^{r_{1}} (I) } \| D^{\alpha} u
\|_{L_{t}^{q_{2}} L_{x}^{r_{2}}(I)}
\end{array}
\label{Eqn:FracLeibn}
\end{equation}
with $ \alpha >0$, $r$, $r_{1}$, $r_{2}$ lying in $[1, \, \infty]$,
$\frac{1}{q}=\frac{1}{q_{1}} + \frac{1}{q_{2}}$ and
$\frac{1}{r}=\frac{1}{r_{1}}+ \frac{1}{r_{2}}$.

The Leibnitz rule for products is

\begin{equation}
\begin{array}{ll}
\| D^{\alpha} (uv) \|_{L_{t}^{q} L_{x}^{r}(I)} & \lesssim  \|
D^{\alpha} u \|_{L_{t}^{q_{1}} L_{x}^{r_{1}}(I)} \| v
\|_{L_{t}^{q_{2}} L_{x}^{r_{2}}(I)} + \| D^{\alpha} u
\|_{L_{t}^{q_{3}} L_{x}^{r_{3}}(I)} \| v  \|_{L_{t}^{q_{4}}
L_{x}^{r_{4}}(I)}
\end{array}
\label{Eqn:Leibn2}
\end{equation}
with $\alpha > 0$, $r$, $r_{1}$, $r_{2}$ lying in $[1, \, \infty]$,
$\frac{1}{q}=\frac{1}{q_{1}} + \frac{1}{q_{2}}$,
$\frac{1}{q}=\frac{1}{q_{3}} + \frac{1}{q_{4}}$, $\frac{1}{r}=
\frac{1}{r_{1}} + \frac{1}{r_{2}}$ and $\frac{1}{r}=\frac{1}{r_{3}}
+ \frac{1}{r_{4}}$.

If $F \in C^{2}$ then we can write

\begin{equation}
\begin{array}{ll}
F(x)-F(y) & = \int_{0}^{1} F^{'}\left( tx+(1-t) y \right).(x-y) \,
dt
\end{array}
\end{equation}
By using (\ref{Eqn:FracLeibn}) and (\ref{Eqn:Leibn2}) the Leibnitz
rule for differences can be formulated as follows

\begin{equation}
\begin{array}{ll}
\| D^{\alpha}( F(u) - F(v) ) \|_{L_{t}^{q} L_{x}^{r}(I)} & \lesssim
\sup_{t \in [0,1]} \| F^{'}(tu+(1-t)v) \|_{L_{t}^{q_{1}}
L_{x}^{r_{1}}(I)} \| D^{\alpha}(u-v) \|_{L_{t}^{q_{2}}
L_{x}^{r_{2}}(I)} \\
&  + \sup_{t \in [0,1]}  \| F^{''}(tu+(1-t)v) \|_{L_{t}^{q^{'}_{1}}
L_{x}^{r^{'}_{1}}(I)}  \left( \|D^{\alpha} u \|_{L_{t}^{q^{'}_{2}}
L_{x}^{r^{'}_{2}}(I)} + \|D^{\alpha} v \|_{L_{t}^{q^{'}_{2}}
L_{x}^{r^{'}_{2}}(I)} \right) \\
&  \| u-v \|_{L_{t}^{q^{'}_{3}} L_{x}^{r^{'}_{3}}(I)}
\end{array}
\label{Eqn:FracLeibnDiff}
\end{equation}
with $\alpha > 0$, $r_{1}$, $r_{2}$, $r^{'}_{1}$, $r^{'}_{2}$,
$r^{'}_{3}$ lying in $[1, \, \infty]$, $\frac{1}{q} =\frac{1}{q_{1}}
+ \frac{1}{q_{2}}$, $\frac{1}{r}=\frac{1}{r_{1}} + \frac{1}{r_{2}}$,
$\frac{1}{q}=\frac{1}{q^{'}_{1}} + \frac{1}{q^{'}_{2}} +
\frac{1}{q^{'}_{3}}$ and $\frac{1}{r}=\frac{1}{r^{'}_{1}} +
\frac{1}{r^{'}_{2}} + \frac{1}{r^{'}_{3}}$.

We shall apply these formulas to several formulas of $F(u)$, and, in
particular, to $F(u):=|u|^{p-1} u g(|u|)$. Notice that, by
(\ref{Eqn:Conditiong2}) and (\ref{Eqn:Conditiong3}), we have
$F^{'}(x) \sim |x|^{p-1} g(|x|)$ and $F^{''}(x) \sim |x|^{p-2}
g(|x|)$. Notice also that, by (\ref{Eqn:Conditiong2}) again, we
have, for $t \in [0,1]$

\begin{equation}
\begin{array}{ll}
g \left( |tx+(1-t)y| \right) & \leq g \left( 2 \max{(|x|,|y|)} \right) \\
& \leq g(\max{(|x|,|y|)} + \log{2} \\
& \lesssim g(|x|) + g(|y|)
\end{array}
\end{equation}
This will allow us to estimate easily $\sup_{t \in [0,1]} \| F^{'} (
tu+(1-t) v ) \|_{L_{t}^{q_{1}} L_{x}^{r_{1}}(I)}$ and $\sup_{t \in
[0,1]} \| F^{''}(tu+(1-t)v) \|_{L_{t}^{q_{1}} L_{x}^{r_{1}}(I)}$.

\vspace{5mm}

Now we explain the main ideas of this paper. We shall prove, in
Section \ref{sec:contLtinf}, that for a large number of $g \, s$, a
special property for the solution of (\ref{Eqn:BarelyWaveEq}) holds

\begin{prop}{\textbf{"control of $S(I)$-norm and control of norm of initial data imply control of $L_{t}^{\infty}
\tilde{H}^{2}(I) \times L_{t}^{\infty} \tilde{H}^{1}(I)$ norm"}} Let
$I \subset I_{max,g}$ and $a \in I$. Assume that $g$ satisfies
(\ref{Eqn:Conditiong2}), (\ref{Eqn:Conditiong3}) and
\footnote{Condition (\ref{Eqn:Conditiong4}) basically says that $g$
grows slowly on average.}

\begin{equation}
\begin{array}{ll}
\int_{1}^{\infty} \frac{1}{y g^{2}(y)} \, dy  =  \infty
\end{array}
\label{Eqn:Conditiong4}
\end{equation}
Let $A \geq 0$ such that $ \| ( u_{0}, u_{1}) \|_{\tilde{H}^{2}
\times \tilde{H}^{1}} \leq A$. Let $u$ be the solution of
(\ref{Eqn:BarelyWaveEq}). There exists a constant $C>0$ such that

\begin{equation}
\begin{array}{ll}
\| ( u ,\partial_{t} u ) \|_{L_{t}^{\infty} \tilde{H}^{2}(I) \times
L_{t}^{\infty} \tilde{H}^{1}(I)} & \leq (2C)^{N} A
\end{array}
\label{Eqn:BoundTildeH2}
\end{equation}
with $N:=N(I)$ such that

\begin{equation}
\begin{array}{ll}
\int_{2CA}^{(2C)^{N} A} \frac{1}{y g^{2}(y)} \, dy
>> \| u \|^{2(p-1)}_{S(I)}
\end{array}
\end{equation}
\label{Prop:controlLtinfty}
\end{prop}
Moreover we shall give a criterion of global well-posedness (proved
in Section \ref{sec:Globcrit})

\begin{prop}{\textbf{"criterion of global well-posedness"}}
Assume that $|I_{max,g}| < \infty$. Assume that $g$ satisfies
(\ref{Eqn:Conditiong2}), (\ref{Eqn:Conditiong3}) and
(\ref{Eqn:Conditiong4}). Then

\begin{equation}
\begin{array}{ll}
 \| u \|_{S(I_{max,g})} & =
\infty
\end{array}
\end{equation}
\label{Prop:CritGlob}
\end{prop}
The first step would be to prove global well-posedness of
(\ref{Eqn:BarelyWaveEq}), with $g_{1}$ a non decreasing function
that is constant for $x$ large (say $x \geq C^{'}_{1}$, with
$C^{'}_{1}$ to be determined). By Proposition \ref{Prop:CritGlob},
it is enough to find an upper bound of the $S([-T,T])$ norm of the
solution $u_{[1]}$ for $T$ arbitrary large. This is indeed possible,
by proving that $g_{1}$ can be considered as a subcritical
perturbation of the nonlinearity. In other words, $g_{1}(|u|)
|u|^{p-1} u$ will play the same role as that of $|u|^{p-1} u \left(
1 - \frac{1}{|u|^{\alpha}} \right)$ for some $\alpha>0$. Once we
have noticed that this comparison is possible, we shall estimate the
relevant norms (and, in particular $ \| u_{[1]} \|_{S([-T,T])}) $ )
by using perturbation theory, Conjecture \ref{Conj:1} and Condition
\ref{Conj:2}, in the same spirit as Zhang \cite{zhang}. We expect to
find a bound of the form

\begin{equation}
\begin{array}{ll}
\| u_{[1]} \|_{S([-T,T])} & \leq C_{3} \left( \| (u_{0},u_{1})
\|_{\tilde{H}^{2} \times \tilde{H}^{1}}, \, T \right)
\end{array}
\label{Eqn:BoundSPrinc}
\end{equation}
with $C_{3}$ increasing as $T$ or $ \| (u_{0},u_{1})
\|_{\tilde{H}^{2} \times \tilde{H}^{1}} $ grows. Notice that if we
restrict $[-T,T]$ to the interval $[-1,1]$ and if the $\tilde{H}^{2}
\times \tilde{H}^{1}$ norm of the initial data $(u_{0},u_{1}) $ is
bounded by one, then we can prove, by (\ref{Eqn:BoundSPrinc}),
(\ref{Eqn:SobEmbed1}) and Proposition \ref{Prop:controlLtinfty},
that the $L_{t}^{\infty} L_{x}^{\infty} ([-T,T])$ norm of the
solution $u_{[1]}$ is bounded by a constant (denoted by $C_{1}$) on
$[-1,1]$. Therefore, if $h$ is a smooth extension of $g_{1}$ outside
$[0,C_{1}]$, and if $u$ is the solution of (\ref{Eqn:BarelyWaveEq})
(with $g:=h$ ) then we expect to prove that $u=u_{[1]}$ on $[-1,1]$
and for data $\|(u_{0},u_{1})\|_{\tilde{H}^{2} \times \tilde{H}^{1}}
\leq 1$. This implies in particular, by (\ref{Eqn:BoundSPrinc}),
that we have a finite upper bound $\| u \|_{S([-1,1])}$. We are not
done yet. There are two problems. First $g_{1}$ does not go to
infinity. Second we only control $\| u \|_{S([-1,1])}$ for data $\|
(u_{0},u_{1}) \|_{\tilde{H}^{2} \times \tilde{H}^{1}} \leq 1$: we
would like to control $\| u \|_{S(\mathbb{R})}$ for arbitrary data.
In order to overcome these difficulties we iterate the procedure
described above. More precisely, given a function $g_{i-1}$ that is
constant for $x \geq C_{i-1}$  and such that $u_{[i-1]}$, solution
of (\ref{Eqn:BarelyWaveEq}) with $g=g_{i-1}$, satisfies  $\|
u_{[i-1]} \|_{S([-(i-1), i-1])} < \infty $, we construct a function
$g_{i}$ that satisfies the following properties

\begin{itemize}

\item it is an extension of $g_{i-1}$ outside $[0,C_{i-1}]$

\item it is increasing and constant (say equal to $i+1$) for $x \geq C^{'}_{i}$, with
$C^{'}_{i}$ to be determined

\end{itemize}
Again, we shall prove that $g_{i}$ may be regarded as a subcritical
perturbation of the nonlinearity $(i+1) |u|^{p-1} u$. This allow us
to control $\| u_{[i]} \|_{S([-i,i])}$, by using perturbation
theory, Conjecture \ref{Conj:1} and Condition \ref{Conj:2}. By using
Proposition \ref{Prop:controlLtinfty} and (\ref{Eqn:SobEmbed1}), we
can find a finite upper bound of $\| u_{[i]} \|_{L_{t}^{\infty}
L_{x}^{\infty} ([-i,i])}$. We assign the value of this upper bound
to $C_{i}$. To conclude the argument we let $\tilde{g}=\lim_{i
\rightarrow \infty} g_{i}$. Given $T > 0$, we can find a $j$ such
that $[-T,T] \subset [-j,j]$ and $\| (u_{0},u_{1}) \|_{\tilde{H}^{2}
\times \tilde{H}^{1}} \leq j$. We prove that $u=u_{[j]}$ on $[-j,j]$
with $u$, solution of (\ref{Eqn:BarelyWaveEq}) with $g:=\tilde{g}$.
Since we have a finite upper bound of $ \| u_{[j]} \|_{S([-j,j])} $,
we also control $\| u \|_{S([-j,j])}$ and  $\| u \|_{S([-T,T])}$.
Theorem \ref{Thm:GlobWave} follows from Proposition
\ref{Prop:CritGlob}.

$\textbf{Acknowledgements}:$ The author would like to thank Terence
Tao for suggesting this problem and for valuable discussions related
to this work.

\section{Proof of Proposition \ref{prop:localexistencebar}}
\label{Sec:Prel}

In this section we prove Proposition \ref{prop:localexistencebar}
for barely $\dot{H}^{s_{p}}(\mathbb{R}^{3})$ supercritical wave
equations (\ref{Eqn:BarelyWaveEq}). The proof is based upon standard
arguments. Here we have chosen to modify an argument in
\cite{kenigmerle1}.

For $\delta$, $T_{l}$, $C$, $M$ to be chosen and such that
(\ref{Eqn:SmallCond}) holds we define

\begin{equation}
\begin{array}{ll}
B_{1} & : = \mathbf{B} \left( \mathcal{C} ( [0,T_{l}], \,
\tilde{H}^{2} ) \cap D^{\frac{1}{2} -s_{p}} W ([0,T_{l}] )
\cap D^{\frac{1}{2}-2} W ([0,T_{l}] ) , 2 C M   \right) \\
B_{2} & := \mathbf{B} \left( S
([0,T_{l}] ), \, 2 \delta \right) \\
B^{'} & : = \mathbf{B} \left( \mathcal{C} ( [0,T_{l}], \,
\tilde{H}^{1} ), 2 CM \right)
\end{array}
\end{equation}
and
\begin{equation}
\begin{array}{ll}
X & : = \left\{ (u,\partial_{t}u): \,  u \in  B_{1} \cap B_{2} , \,
\partial_{t} u \in B^{'}  \right\}
\end{array}
\end{equation}
Let

\begin{equation}
\begin{array}{ll}
\Psi (u, \, \partial_{t} u ) & := \left(
\begin{array}{l}
\cos{(tD)} u_{0} + \frac{\sin{(tD)}}{D} u_{1} - \int_{0}^{t}
\frac{\sin{(t-t^{'})D}}{D} \left( |u(t^{'})|^{p-1} u(t^{'}) g(|u(t^{'})|) \right) \, d t^{'}  \\
-D \sin{(tD)} u_{0} + \cos{(tD)} u_{1} - \int_{0}^{t} \
\cos{(t-t^{'}) D}  \left( |u(t^{'})|^{p-1} u(t^{'}) g(|u(t^{'})|)
\right) \, d t^{'}
\end{array}
\right)
\end{array}
\end{equation}
Then

\begin{itemize}

\item $\Psi$ maps $X$ to $X$. Indeed, by (\ref{Eqn:SmallCond}),
(\ref{Eqn:Strichartz}) and the fractional Leibnitz rule
(\ref{Eqn:FracLeibn}) applied to $\alpha \in \left\{ s_{p} -
\frac{1}{2}, \, 2 -\frac{1}{2} \right\}$ and $F(u):= |u|^{p-1} u
g(|u|)$ and by applying the multipliers $D^{2-\frac{1}{2}}$ and
$D^{s_{p}-\frac{1}{2}}$ to the Strichartz estimates with
$m=\frac{1}{2}$ we have

\begin{equation}
\begin{array}{ll}
Q \left( [0, \, T_{l}] \right) & \lesssim \left\| ( u_{0}, \, u_{1})
\right\|_{\tilde{H}^{2}(\mathbb{R}^{3}) \times
\tilde{H}^{1}(\mathbb{R}^{3})} + \| D^{s_{p} - \frac{1}{2}} \left(
|u|^{p-1} u g ( |u| ) \right) \|_{\tilde{W}([0,T_{l}])}    \\
& + \| D^{2 - \frac{1}{2}}
\left( |u|^{p-1} u g ( |u| ) \right) \|_{\tilde{W}([0,T_{l}])} \\ & \\
& \leq CM + C \left( \|D^{s_{p}-\frac{1}{2}} u \|_{W([0,T_{l}])} +
\| D^{2-\frac{1}{2}} u \|_{W([0,T_{l}])} \right) \| u
\|^{p-1}_{S([0,T_{l}])}  g ( \| u
\|_{L_{t}^{\infty}L_{x}^{\infty}([0,T_{l}])} )
\\ & \\ & \leq
CM + (2 \delta)^{p-1} C (2CM)  g(2CM)
\end{array}
\end{equation}
for some $C>0$ and

\begin{equation}
\begin{array}{ll}
\| u \|_{S ( [0,T_{l}])} - \delta & \lesssim \|
D^{s_{p}-\frac{1}{2}} \left( |u|^{p-1} u g(|u|) \right)
\|_{\tilde{W}
([0,T_{l}] )} \\
& \lesssim \| u \|^{p-1}_{S ([0, T_{l}]) } \| D^{s_{p} -\frac{1}{2}}
u \|_{W ([0,T_{l}])}  g ( \| u \|_{L_{t}^{\infty} L_{x}^{\infty} (
[0,T_{l}]) })  \\ &  \\ & \lesssim   (2 \delta)^{p-1} (2CM) g(2CM)
\end{array}
\end{equation}
Choosing $\delta=\delta(M)
> 0$ small enough we see that $\Psi(X) \subset X$.

\item $\Psi$ is a contraction. Indeed we get

\begin{equation}
\begin{array}{l}
\| \Psi (u) - \Psi(v)  \|_{X} \\
 \lesssim \left\|  D^{s_{p}
-\frac{1}{2}} (|u|^{p-1}u g(|u|) - |v|^{p-1} v  g(|v|) ) \right
\|_{\tilde{W}([0,T_{l}])} + \left\|  D^{2 -\frac{1}{2}} ( |u|^{p-1}u
g(|u|) - |v|^{p-1} v
g(|v|) ) \right\|_{\tilde{W}([0,T_{l}])} \\ \\
\lesssim g ( \| u \|_{L_{t}^{\infty} L_{x}^{\infty}([0,T_{l}])}) \times \\
\left(
\begin{array}{l}
\left( \| u \|^{p-1}_{S([0,T_{l}])} + \| v \|^{p-1}_{S([0,T_{l}])}
\right) \left( \| D^{s_{p}- \frac{1}{2}} (u -v) \|_{W([0,T_{l}])}
+ \| D^{2-\frac{1}{2}} (u-v) \| _{W([0,T_{l}])} \right) \\
+ \left( \| u \|^{p-2}_{S([0,T_{l}])} + \| v \|^{p-2}_{S([0,T_{l}])}
\right) \left(
\begin{array}{l}
\| D^{s_{p}-\frac{1}{2}} u \|_{W([0,T_{l}])} + \| D^{2-\frac{1}{2}}
u \|_{W([0,T_{l}])} + \| D^{s_{p}-\frac{1}{2}} v \|_{W([0,T_{l}])} \\
 + \| D^{2-\frac{1}{2}} v \|_{W([0,T_{l}])}
\end{array}
\right)
\\
\| u- v \|_{S([0,T_{l}])}
\end{array}
\right) \\ \\
 \lesssim \left( g(2CM) (2 \delta)^{p-1} + (2 \delta)^{p-2} (2CM)
\right) \| u - v \|_{X}
\end{array}
\label{Eqn:DiffPsi}
\end{equation}
In the above computations, we applied the Leibnitz rule for
differences to $\alpha \in \{ s_{p}-\frac{1}{2}, 2-\frac{1}{2} \}$
and $F(u):=|u|^{p-1} u g(|u|)$. Therefore, if $\delta=\delta(M)
> 0$ is small enough, then $\Psi$ is a contraction.

\end{itemize}

\section{Proof of Proposition \ref{Prop:controlLtinfty}}
\label{sec:contLtinf}

In this section we prove Proposition \ref{Prop:controlLtinfty}.

It is enough to prove that $Q(I) < \infty $. Without loss of
generality we can assume that $A >>1$. Then we divide $I$ into
subintervals $(I_{i})_{1 \leq i \leq N}$ such that

\begin{equation}
\begin{array}{ll}
\| u \|_{S(I_{i})} & = \frac{\eta}{g^{\frac{1}{p-1}} \left( (2
C)^{N} A \right)}
\end{array}
\end{equation}
for some $C \gtrsim 1$ and $\eta > 0$ constants to be chosen later,
except maybe the last one. Notice that such a partition always
exists since by (\ref{Eqn:Conditiong4}) we get for $N:=N(I)$ large
enough

\begin{equation}
\begin{array}{ll}
\sum_{i=1}^{N} \frac{1}{ g^{2} ( (2 C)^{i} A) } &
\geq \int_{1}^{N} \frac{1}{ g^{2}((2C)^{x} A )} \, dx  \\
& \gtrsim \int_{2CA}^{(2C)^{N} A} \frac{1}{y g^{2}(y)} \,
dy \\
& >> \| u \|^{2(p-1)}_{S(I)}
\end{array}
\end{equation}
We get, by a similar token used in Section \ref{Sec:Prel}

\begin{equation}
\begin{array}{ll}
Q(I_{1},u) & \lesssim \| (u_{0},u_{1})
\|_{\tilde{H}^{2}(\mathbb{R}^{3}) \times
\tilde{H}^{1}(\mathbb{R}^{3})} + \left\| D^{s_{p}- \frac{1}{2}}
\left( |u|^{p-1} u g(|u|) \right) \right\|_{\tilde{W}(I_{1})} \\
& +
\left\| D^{2- \frac{1}{2}} \left( |u|^{p-1} u g(|u|) \right) \right\|_{\tilde{W}(I_{1})} \\
& \lesssim A + \left( \| D^{s_{p} -\frac{1}{2}} u \|_{W(I_{1})} + \|
 D^{2-\frac{1}{2}} u \|_{W(I_{1})} \right) \| u
\|^{p-1}_{S(I_{1})} g ( \| u \|_{ L_{t}^{\infty} L_{x}^{\infty}(I_{1})} )  \\
& \lesssim A + \| u \|_{S(I_{1})}^{p-1} Q(I_{1},u) g( Q(I_{1},u) )
\end{array}
\label{Eqn:IneqQIi}
\end{equation}
We choose $C$ to be equal to the constant determined by $\lesssim$
in (\ref{Eqn:IneqQIi}). Without loss of generality we can assume
that $C>1$. By a continuity argument, iteration on $i$, we get, for
$\eta << 1$, (\ref{Eqn:BoundTildeH2}).

\section{Proof of Proposition \ref{Prop:CritGlob}}
\label{sec:Globcrit}

In this section we prove Proposition \ref{Prop:CritGlob}.

We argue as follows: by time reversal symmetry it is enough to prove
that $T_{+,g} < \infty$. If $ \| u \|_{S(I_{max,g})} < \infty $ then
we have $ Q([0,T_{+,g}],u) < \infty $: this follows by slightly
adapting the proof of Proposition \ref{Prop:controlLtinfty}.
Consequently, by the dominated convergence theorem, there would
exists a sequence $t_{n} \rightarrow T_{+,g}$ such that $\| u \|_{
S( [t_{n}, T_{+,g}] ) } << \delta $ and $ \| D^{ s_{p} -
\frac{1}{2}} u \|_{W([t_{n},T_{+,g}])} << \delta $ if $n$ is large
enough with $\delta$ defined in Proposition
\ref{prop:localexistencebar}. But, by (\ref{Eqn:DuhForm}) and
(\ref{Eqn:Strichartz})

\begin{equation}
\begin{array}{l}
\| \cos{ \left( (t-t_{n}) D \right)} u(t_{n}) + \frac{\sin{
(t-t_{n}) D }}{ D} u_{1} \|_{S( [t_{n},T_{+,g}] )} \\
 \lesssim \| u
\|_{S( [t_{n},T_{+,g}] )}   + \| u \|_{S([t_{n}, T_{+,g}])}^{p-1} \|
D^{s_{p} -\frac{1}{2}} u \|_{W([t_{n}, T_{+,g}])}  g (
Q([0,T_{+,g},u]) ) \\   << \delta
\end{array}
\end{equation}
and consequently, by continuity, there would exist $\tilde{T} >
T_{+,g}$ such that

\begin{equation}
\begin{array}{ll}
\| \cos{ \left( (t-t_{n}) D \right)} u(t_{n}) + \frac{\sin{(t-t_{n})
D}}{ D} \partial_{t} u(t_{n}) \|_{S( [t_{n}, \tilde{T}] ) } & \leq
\delta
\end{array}
\end{equation}
, which would contradict the definition of $T_{+,g}$.

\begin{rem}
Notice that if we have the stronger bound $ \| u \|_{S(I_{max,g})}
\leq C $ with $C:=C \left( \| ( u_{0},u_{1}) \|_{\tilde{H}^{2}
\times \tilde{H}^{1}} \right) < \infty$, then not only $I_{max,g} =
(- \infty, + \infty)$ but also $u$ scatters as $t \rightarrow \pm
\infty$. Indeed by Proposition \ref{Prop:CritGlob},
$I_{max,g}=\mathbb{R}$. Then by time reversal symmetry it is enough
to assume that $t \rightarrow \infty$. Let $v(t):= \left( u(t),
\partial_{t} u(t) \right)$.  We are looking for $v_{+,0}:= \left(
u_{+,0}, \, u_{+,1} \right)$ such that

\begin{equation}
\begin{array}{ll}
\left\| v(t) - K (t) v_{+,0} \right\|_{\tilde{H}^{2} \times
\tilde{H}^{1}} & \rightarrow 0
\end{array}
\end{equation}
as $t \rightarrow  \infty$. Here

\begin{equation}
\begin{array}{ll}
K(t) & : = \left(
\begin{array}{ll}
\cos{tD} &  \frac{\sin{tD}}{D} \\
-D \sin{tD} & \cos{tD}
\end{array}
\right)
\end{array}
\end{equation}
We have

\begin{equation}
\begin{array}{ll}
K^{-1} (t) & = \left(
\begin{array}{ll}
\cos{(tD)} & -\frac{\sin{(tD)}}{D} \\
D \sin{(tD)} & \cos{(tD)}
\end{array}
\right)
\end{array}
\end{equation}
Notice that $K^{-1}(t)$ and $K(t)$ are bounded in $ \tilde{H}^{2}
\times \tilde{H}^{1} $. Therefore it is enough to prove that
$K^{-1}(t) v(t) $ has a limit as $t \rightarrow  \infty$. But since
$K^{-1}(t) v(t) = (u_{0},u_{1}) - K^{-1}(t) \left( u_{nl}(t), \,
\partial_{t} u_{nl}(t) \right)$ (with $u_{nl}$ denoting the
nonlinear part of the solution (\ref{Eqn:Duhamelg}) \footnote{i.e
\begin{equation}
\begin{array}{ll}
u_{nl}(t)  & :=- \int_{0}^{t} \frac{\sin{(t-t^{'}) D}}{D} \left(
|u(t^{'})|^{p-1} u(t^{'}) g(|u(t^{'})|)  \right) \, dt^{'}
\end{array}
\nonumber
\end{equation}
}), then it suffices to prove that $K^{-1}(t) \left( u_{nl}(t),
\partial_{t} u_{nl}(t) \right)$ has a limit. But

\begin{equation}
\begin{array}{l}
\| K^{-1} (t_{1})  u_{nl}(t_{1}) - K^{-1}(t_{2}) u_{nl}(t_{2})
\|_{\tilde{H}^{2} \times
\tilde{H}^{1}} \\ \\
 \lesssim \left\| ( u_{nl},
\partial_{t} u_{nl} ) \right\|_{ L_{t}^{\infty} \tilde{H}^{2} ([t_{1},t_{2}]) \times
L_{t}^{\infty} \tilde{H}^{1} ([t_{1},t_{2}]) } \\ \\
 \lesssim  \left( \| D^{s_{p} -\frac{1}{2}}( |u|^{p-1} u g(|u|)  )
\|_{\tilde{W}([t_{1},t_{2}])} + \| D^{2 -\frac{1}{2}}( |u|^{p-1} u
g(|u|)  ) \|_{\tilde{W}([t_{1},t_{2}])} \right) \\ \\
\lesssim \left(
\| D^{s_{p} -\frac{1}{2}} u \|_{W([t_{1},t_{2}])} + \|D^{ 2
-\frac{1}{2}} u \|_{W([t_{1},t_{2}])} \right) \| u
\|^{p-1}_{S([t_{1},t_{2}])}
g( \| u \|_{L_{t}^{\infty}  L_{x}^{\infty} (\mathbb{R}) } )   \\
\end{array}
\end{equation}
It remains to prove that $Q (\mathbb{R}) < \infty$  in order to
conclude that the Cauchy criterion is satisfied, which would imply
scattering. This follows from  $\| u \|_{S(\mathbb{R})} < \infty$
and a slight modification of the proof of Proposition
\ref{Prop:controlLtinfty}.
\end{rem}

\section{Construction of the function $\tilde{g}$}

In this section we prove Theorem \ref{Thm:GlobWave}. Let

\begin{equation}
\begin{array}{ll}
Up \, (i) : = \left\{ (T, \, (u_{0}, u_{1})  ): \, 0 \leq  T \leq i,
\, \| (u_{0},u_{1}) \|_{\tilde{H}^{2} \times \tilde{H}^{1}} \leq i
\right\}
\end{array}
\end{equation}
As $i$ ranges over $\{1,2,... \}$ we construct, for each set $ Up \,
(i) $, a function $g_{i}$ satisfying (\ref{Eqn:Conditiong2}) and
(\ref{Eqn:Conditiong3}). Moreover it is constant for large values of
$|x|$. The function $g_{i+1}$ depends on $g_{i}$; the construction
of $g_{i}$ is made by induction on $i$. More precisely we will prove
that the following lemma

\begin{lem}

Let $A >>1$. Then there exist two sequences of numbers $\{ C_{i}
\}_{i \geq 0}$, $ \{ C^{'}_{i} \}_{i \geq 0}$ and a sequence of
functions $ \{ g_{i} \}_{i \geq 0}$ such that for all $\left( T, \,
(u_{0},u_{1}) \right) \in Up \, (i)$

\begin{itemize}

\item $g_{0}:=1$, $C_{0}:=0$ and $C^{'}_{0}=0$

\item $\{ C_{i} \}_{i \geq 0}$ and  $\{ C^{'}_{i} \}_{i \geq 0}$ are positive, non decreasing. Moreover they satisfy

\begin{equation}
\begin{array}{l}
A C_{i-1} < C^{'}_{i} < A C_{i}
\end{array}
\label{Eqn:Ci1}
\end{equation}
for $i \geq 1$ and

\begin{equation}
\begin{array}{ll}
C_{i} &  \geq  i
\end{array}
\label{Eqn:Ci2}
\end{equation}

\item $g_{i}$ is smooth, non decreasing; it satisfies (\ref{Eqn:Conditiong2}),
(\ref{Eqn:Conditiong3}),

\begin{equation}
\begin{array}{ll}
\int_{1}^{C^{'}_{i}} \frac{1}{t  g_{i}^{2}(t)} \, dt &
\rightarrow_{i \rightarrow \infty}  \infty
\end{array}
\end{equation}
and

\begin{equation}
\left\{
\begin{array}{ll}
g_{i}(|x|) & = g_{i-1}(|x|), \, |x| \leq A C_{i-1} \\
g_{i}(|x|) & = i+1, \, |x| \geq C^{'}_{i}
\end{array}
\right. \label{Eqn:Constructiongi}
\end{equation}

\item the solution $u_{[i]}$ of the following wave equation

\begin{equation}
\left\{
\begin{array}{ll}
\partial_{tt} u_{[i]} - \triangle u_{[i]} & = - |u_{[i]}|^{p-1}
u_{[i]} g_{i} \left( |u_{[i]}| \right) \\
u_{[i]}(0) & = u_{0} \in \tilde{H}^{2} \\
\partial_{t} u_{[i]}(0) & = u_{1} \in \tilde{H}^{1}
\end{array}
\right.
\label{Eqn:Wavegi}
\end{equation}
satisfies

\begin{equation}
\begin{array}{ll}
 \max{ \left( \| u_{[i]} \|_{S([-i,i])} , \left\| ( u_{[i]}, \partial_{t} u_{[i]}) \right\|_{L_{t}^{\infty} \tilde{H}^{2} ([-T,T])
\times L_{t}^{\infty} \tilde{H}^{1}([-T,T]) } \right) } & \leq C_{i}
\end{array}
\label{Eqn:Boundpropi}
\end{equation}
\end{itemize}
\label{lem:gi}
\end{lem}

This lemma will be proved in the next subsection. Assuming that is
true let $\tilde{g} = \lim_{i \rightarrow \infty} g_{i}$. Then
clearly $\tilde{g}$ is smooth; it satisfies (\ref{Eqn:Conditiong2})
and (\ref{Eqn:Conditiong3}). It also goes to infinity. Moreover let
$u$ be the solution of (\ref{Eqn:BarelyWaveEq}) with $g:=\tilde{g}$.
We want to prove that the solution $u$ exists for all time. Let
$T_{0} \geq 0$ be a fixed time. Let $j:=j(T_{0},\| u_{0}
\|_{\tilde{H}^{2}}, \| u_{1} \|_{\tilde{H}^{1}})>0$ be the smallest
positive integer such that $[-T_{0},T_{0}] \subset [-j,j]$ and  $\|
(u_{0},u_{1} ) \|_{\tilde{H}^{2} \times \tilde{H}^{1}} \leq j$. We
claim that $\|( u, \partial_{t} u) \|_{L_{t}^{\infty}
\tilde{H}^{2}([-T_{0},T_{0}]) \times L_{t}^{\infty}
\tilde{H}^{1}([-T_{0},T_{0}]) } \leq C_{j}$ and $\| u
\|_{S([-T_{0},T_{0}])} \leq C_{j}$. Indeed let

\begin{equation}
\begin{array}{ll}
F_{j} & : = \left\{ t \in [0,j], \, \|(u, \partial_{t} u)
\|_{L_{t}^{\infty} \tilde{H}^{2}([-t,t]) \times L_{t}^{\infty}
\tilde{H}^{1}([-t,t]) } \leq C_{j} \, \mathrm{and} \, \| u
\|_{S([-t,t])} \leq C_{j} \right\}
\end{array}
\end{equation}
Then

\begin{itemize}

\item $F_{j} \neq \emptyset$: $0 \in  F_{j}$, by (\ref{Eqn:Ci2}).

\item $F_{j}$ is closed. Indeed let $\tilde{t} \in  \overline{F_{j}}$. Then there exist a sequence $(t_{n})_{n \geq 1}$ such that
$t_{n} \in [0,j]$ such that $t_{n} \rightarrow \tilde{t}$, $ \| u
\|_{S([-t_{n},t_{n}])} \leq C_{j}$ and \\ $\| (  u,\partial_{t} u )
\|_{L_{t}^{\infty} \tilde{H}^{2} ([-t_{n},t_{n}]) \times
L_{t}^{\infty} \tilde{H}^{1} ([-t_{n},t_{n}]) }   \leq C_{j}$. It is
enough to prove that $ \| u \|_{S([-\tilde{t}, \tilde{t}])} <
\infty$ and then apply dominated convergence. There are two
subcases:

\begin{itemize}

\item $\card \{t_{n}, \, t_{n} \leq \tilde{t} \} < \infty $:
in this case, there exists $n_{0}$ large enough such that $t_{n}
\geq \tilde{t}$ for $n \geq n_{0}$ and

\begin{equation}
\begin{array}{ll}
\| u \|_{S([-\tilde{t}, \tilde{t}])} & \leq \| u
\|_{S([-t_{n},t_{n}])} \\
& < \infty
\end{array}
\end{equation}

\item  $\card \{t_{n}, \, t_{n} \leq \tilde{t} \} = \infty $:
even if it means passing to a subsequence, we can still assume that
$t_{n} \leq \tilde{t}$. Let $n_{0} \geq 1$ be fixed. Then, by the
inequality
\begin{equation}
\begin{array}{ll}
\| \cos{(t-t_{n_{0}})D} u(t_{n_{0}}) + \frac{\sin{(t-t_{n_{0}})
D}}{D} \partial_{t} u(t_{n_{0}}) \|_{S([t_{n_{0}}, \tilde{t}])} &
\lesssim \left\| ( u(t_{n_{0}}), \partial_{t} u(t_{n_{0}}) )
\right\|_{\tilde{H}^{2} \times \tilde{H}^{1}} \\
& \lesssim C_{j}
\end{array}
\end{equation}
we conclude from the dominated convergence theorem that there exists
$n_{1}:= n_{1}(n_{0})$ large enough such that 
\begin{equation}
\begin{array}{ll}
\|\cos{(t-t_{n_{0}})D} u(t_{n_{0}}) + \frac{\sin{(t-t_{n_{0}})
D}}{D} \partial_{t} u(t_{n_{0}}) \|_{S([t_{n_{1}}, \tilde{t} ])} &
\leq \delta
\end{array}
\label{Eqn:ControlLin}
\end{equation}
with $\delta:=\delta(C_{j})$ defined in Proposition
\ref{prop:localexistencebar}. Therefore, by Proposition
\ref{prop:localexistencebar}, we have $\| u \|_{S([t_{n_{1}},
\tilde{t}])} < \infty$. By a similar token, $\| u \|_{S([-\tilde{t},
-t_{n_{1}}])} < \infty$. Combining these inequalities with $\| u
\|_{S([-t_{n_{1}}, t_{n_{1}}])} \leq C_{j} $, we eventually get $\|
u \|_{S([-\tilde{t},\tilde{t}])} < \infty$.
\end{itemize}

\item $F_{j}$ is open. Indeed let $\bar{t} \in F_{j}$. By Proposition \ref{prop:localexistencebar} there exists
$\alpha  > 0$ such that if $ t \in (\bar{t} - \alpha, \bar{t} +
\alpha) \cap [0,j] $  then $ [-t,t] \subset I_{max,\tilde{g}} $ and
$ \| u \|_{L_{t}^{\infty} L_{x}^{\infty} ([-t,t])} \lesssim \| u
\|_{L_{t}^{\infty} \tilde{H}^{2}([-t,t])} \lesssim   C_{j} $. Notice
also that, by (\ref{Eqn:Boundpropi}),  $ [-t,t] \subset I_{max,
g_{j}} $. In view of these remarks, we conclude, after slightly
adapting the proof of Proposition \ref{Prop:controlLtinfty}, that $
Q( [-t,t], u ) \lesssim_{j} 1 $ and $Q([-t,t], u_{[j]}) \lesssim_{j}
1 $. We divide $[-t, \, t] $ into a finite number of subintervals
$(I_{i})_{i \leq k}=([a_{i},b_{i}])_{1 \leq i \leq k}$ that satisfy,
for $\eta <<1$ to be defined later, the following properties hold

\begin{enumerate}

\item $1 \leq i \leq k$:  $ \| u_{[j]} \|_{S(I_{i})} \leq \eta$,  $ \| u \|_{S(I_{i})} \leq
\eta$, $ \| D^{s_{p}-\frac{1}{2}} u_{[j]} \|_{W(I_{i})} \leq \eta$,
$ \| D^{2-\frac{1}{2}}  u  \|_{W(I_{i})} \leq \eta $ and $ \|
D^{2-\frac{1}{2}} u_{[j]} \|_{W(I_{i})} \leq \eta$

\item $1 \leq i < k$: $ \| u_{[j]}  \|_{S(I_{i})} = \eta $ or $ \| u \|_{S(I_{i})} = \eta
$ or  $ \| D^{s_{p}-\frac{1}{2}} u_{[j]} \|_{W(I_{i})} = \eta$ or $
\| D^{2-\frac{1}{2}}  u  \|_{W(I_{i})} = \eta $ or $ \|
D^{2-\frac{1}{2}} u_{[j]} \|_{W(I_{i})} = \eta$
\end{enumerate}
Notice that, by (\ref{Eqn:Conditiong2}), we have
\begin{equation}
\begin{array}{ll}
\| g_{j}(|u|) - g_{j}(|u_{[j]}|) \|_{L_{t}^{\infty} L_{x}^{\infty}
(I_{i})} & \lesssim \| u-u_{[j]} \|_{L_{t}^{\infty} L_{x}^{\infty}
(I_{i})} \\
& \lesssim \| u-u_{[j]} \|_{L_{t}^{\infty} \tilde{H}^{2} (I_{i})}
\end{array}
\label{Eqn:PtwiseDiffg}
\end{equation}
We consider $w=u-u_{[j]}$. By applying the Leibnitz rules
(\ref{Eqn:FracLeibn}), (\ref{Eqn:FracLeibnDiff}) and
(\ref{Eqn:Leibn2}), and by (\ref{Eqn:PtwiseDiffg}) we have

\begin{equation}
\begin{array}{ll}
Q(I_{1},w) & \lesssim  \left\| D^{s_{p} -\frac{1}{2}} ( |u|^{p-1} u
(\tilde{g}-g_{j})(|u|) )  \right\|_{\tilde{W}(I_{1})} + \left\| D^{2
-\frac{1}{2}} ( |u|^{p-1} u (\tilde{g}-g_{j})(|u|) )
\right\|_{\tilde{W}(I_{1})}  \\
& +  \left\| D^{s_{p} -\frac{1}{2}} ( |u|^{p-1} u - |u_{[j]}|^{p-1}
u_{[j]} ) g_{j}(|u|) \right\|_{\tilde{W}(I_{1})} + \\
& \left\| D^{2 -\frac{1}{2}} ( |u|^{p-1} u
- |u_{[j]}|^{p-1} u_{[j]} ) g_{j}(|u|) \right\|_{\tilde{W}(I_{1})} \\
& + \left\| D^{s_{p} -\frac{1}{2}} ( |u_{[j]}|^{p-1}u_{[j]} (
g_{j}(|u|)-g_{j}(|u_{[j]}|) ) ) \right\|_{\tilde{W}(I_{1})} + \\
& \left\| D^{2 -\frac{1}{2}} ( |u_{[j]}|^{p-1}u_{[j]} ( g_{j}(|u|)-g_{j}(|u_{[j]}|) )) \right\|_{\tilde{W}(I_{1})}    \\
& \lesssim (\tilde{g}-g_{j}) ( \| u \|_{L_{t}^{\infty}
\tilde{H}^{2}(I_{1})} ) \left( \| D^{s_{p} -\frac{1}{2}} u
\|_{W(I_{1})} +  \| D^{2 -\frac{1}{2}} u \|_{W(I_{1})}  \right) \| u
\|^{p-1}_{S(I_{1})}
\\
& + g_{j}( \| u \|_{L_{t}^{\infty} \tilde{H}^{2}(I_{1}) }) \left[
\begin{array}{l}
( \| u_{[j]} \|^{p-1}_{S(I_{1})} + \| u \|^{p-1}_{S(I_{1})} ) (  \|
D^{s_{p}-\frac{1}{2}} w \|_{W(I_{1})}   + \| D^{2-\frac{1}{2}} w \|_{W(I_{1})}  )  \\
+ ( \| u_{[j]} \|^{p-2}_{S(I_{1})} + \| u \|^{p-2}_{S(I_{1})} ) \\
\left(
\begin{array}{l}
\| D^{s_{p}-\frac{1}{2}} u \|_{W(I_{1})} + \| D^{2-\frac{1}{2}} u
\|_{W(I_{1})}  \\
+  \| D^{s_{p}-\frac{1}{2}} u_{[j]} \|_{W(I_{1})}  + \|
D^{2-\frac{1}{2}} u_{[j]} \|_{W(I_{1})}
\end{array}
\right) \| w \|_{S(I_{1})}
\end{array}
\right] \\
& + \| g^{'}_{j}( |u| ) \|_{L_{t}^{\infty} L_{x}^{\infty}(I_{1})}
\left( \| D^{s_{p} -\frac{1}{2}} u \|_{W(I_{1})} + \|
D^{2-\frac{1}{2}} u \|_{W(I_{1})} \right) ( \| u \|^{p-2}_{S(I_{1})}
+
\| u_{[j]} \|^{p-2}_{S(I_{1})} ) \| w \|_{S(I_{1})} \\
& + \left\| g_{j}( | u | ) - g_{j} ( |u_{[j]}|)
\right\|_{L_{t}^{\infty} L_{x}^{\infty}(I_{1})} \| D^{s_{p}
-\frac{1}{2}} u_{[j]} \|_{W(I_{1})} \| u_{[j]} \|^{p-1}_{S(I_{1})} \\
& + \| u_{[j]} \|^{p-1}_{S(I_{1})} \| u_{[j]} \|_{L_{t}^{\infty}
\tilde{H}^{2}(I_{1})} \left(
\begin{array}{l}
\| w \|_{L_{t}^{\infty} \tilde{H}^{2}(I_{1})}  \left( \|
D^{s_{p}-\frac{1}{2}} u \|_{W(I_{1})} + \| D^{s_{p}-\frac{1}{2}}
u_{[j]} \|_{W(I_{1})}
\right) \\
+ \|D^{s_{p}-\frac{1}{2}} w \|_{W(I_{1})}
\end{array} \right)
\\
 & \lesssim g_{j}(C_{j}) \eta^{p-1} Q(I_{1},w) + \eta^{p-1} Q(I_{1},w)
 +\eta^{p} Q(I_{1},w) + C_{j} \eta^{p-1} \left (\eta Q(I_{1},w)  + Q(I_{1},w) \right)
\end{array}
\end{equation}
since, by choosing $A$ large enough and by construction of
$\tilde{g}$ we have

\begin{equation}
\begin{array}{ll}
(\tilde{g}-g_{j})( \| u \| _{L_{t}^{\infty} \tilde{H}^{2}(I_{1})})
& =0
\end{array}
\end{equation}
Therefore we conclude from a continuity argument that $Q(I_{i},w)=0$
and $u=u_{[1]}$ on $I_{1}$. In particular $u(b_{1})=u_{[j]}(b_{1})$.
By iteration on $i$, it is not difficult to see that $u=u_{[1]}$ on
$[-t,t]$. Therefore  $(\bar{t} - \alpha, \bar{t} + \alpha) \cap
[0,j] \subset F_{j}$, by (\ref{Eqn:Boundpropi}).
\end{itemize}

Therefore, $F_{j}=[-j,j]$ and we have $ \| u \|_{S([-T_{0},T_{0}])}
\leq C_{j}$. This proves global well-posedness. Moreover, since $j$
depends in $T_{0}$ and $\|(u_{0},u_{1}) \|_{\tilde{H}^{2} \times
\tilde{H}^{1}}$ , we get (\ref{Eqn:BounduLinftyf}).

\subsection{Proof of Lemma \ref{lem:gi}}

In this subsection we prove Lemma \ref{lem:gi}. We establish
\textit{a priori} bounds.

\begin{itemize}

\item \textit{Step 1}: We construct $g_{1}$. $g_{1}$ is basically a
nonnegative function that increases and is equal to two for $x$
large. Recall that $ [-T,T] \subset [-1,1]$ and $ \|  ( u_{0}, \,
u_{1}) \|_{\tilde{H}^{2} \times \tilde{H}^{1}} \leq 1$. Let $I
\subset [-T,T]$.

Observe that the point $(\infty - ,3+): = \left( \frac{3+
\epsilon}{\epsilon}, 3+ \epsilon  \right)$ with $\epsilon << 1$ is
$\frac{1}{2}$- wave admissible.

\begin{rem}
We would like to chop an interval $I$ such that $\| .
\|_{L_{t}^{\infty} L_{x}^{3}(I)} < \infty$ into subintervals $I_{j}$
such that $\| . \|_{L_{t}^{\infty} L_{x}^{3}(I_{j})}$ is as small as
wanted. Unfortunately this is impossible because the
$L_{t}^{\infty}$ norm is pathologic. Instead we will apply this
process to $\| . \|_{L_{t}^{\infty -} L_{x}^{3+}}$. This creates
slight variations almost everywhere in the process of the
construction of $g_{i}$. Details with respect to these slight
perturbations have been omitted for the sake of readability: they
are left to the reader, who should ignore the $+$ and $-$ signs at
the first reading.
\end{rem}

We define

\begin{equation}
\begin{array}{ll}
X(I) & := D^{\frac{1}{2} -s_{p}} L_{t}^{\infty -} L_{x}^{3 +} (I)
\cap D^{\frac{1}{2} -s_{p}}W(I) \cap S(I) \cap L_{t}^{\infty}
\dot{H}^{s_{p}}(I) \times L_{t}^{\infty} \dot{H}^{s_{p}-1}(I)
\end{array}
\end{equation}
Let $g_{1}$ be a smooth function, defined on the set of nonnegative
real numbers nondecreasing and such that $h_{1}:= g_{1} - 2$
satisfies the following properties: $h_{1}(0)=-1$, h is non
decreasing and $h_{1}(x)=0$ if $|x| \geq 1$. It is not difficult to
see that (\ref{Eqn:Conditiong2}) and (\ref{Eqn:Conditiong3}) are
satisfied.

Observe that

\begin{equation}
\begin{array}{ll}
|h_{1}(x)| & \lesssim \frac{1}{|x|^{\frac{p-1}{2}-}}
\end{array}
\label{Eqn:hpne}
\end{equation}
and

\begin{equation}
\begin{array}{ll}
|h_{1}^{'}(x)| & \lesssim \frac{1}{|x|^{\frac{p+1}{2}-}}
\end{array}
\label{Eqn:hprimeonecd}
\end{equation}
Let $u_{[1]}$ and $v_{[1]}$ be solutions of the equations
\begin{equation}
\left\{
\begin{array}{ll}
\partial_{tt} u_{[1]} - \triangle u_{[1]} & = - \left| u_{[1]} \right|^{p-1} u_{[1]} g_{1} \left( \left| u_{[1]} \right| \right)  \\
u_{[1]}(0) & =u_{0} \in \tilde{H}^{2} \\
\partial_{t} u_{[1]}(0) & = u_{1} \in \tilde{H}^{1}
 \end{array}
\right.
\end{equation}
and

\begin{equation}
\left\{
\begin{array}{ll}
\partial_{tt} v_{[1]} - \triangle v_{[1]} & = - 2 \left| v_{[1]} \right|^{p-1} v_{[1]} \\
v_{[1]}(0) & =u_{0}  \\
\partial_{t} v_{[1]}(0) & = u_{1}
 \end{array}
\right.
\end{equation}
There are several substeps

\begin{itemize}

\item \textit{Substep1}: We claim that $\| v_{[1]} \|_{X(\mathbb{R})}
< \infty$. Indeed, since we assumed that Conjecture \ref{Conj:1} is
true, we first divide $\mathbb{R}$ into subintervals
$(I_{j}=[t_{j},t_{j+1}])_{1 \leq j \leq l}$ such that $\| v_{[1]}
\|_{S(I_{j})} = \eta$ and $\| v_{[1]} \|_{S(I_{l})} \leq \eta$ with
$\eta << 1$. Then we have

\begin{equation}
\begin{array}{ll}
\| v_{[1]} \|_{X(I_{j+1})} & \lesssim \| (
v_{[1]}(t_{j}),\partial_{t} v_{[1]}(t_{j}) )
\|_{\dot{H}^{s_{p}}(\mathbb{R}^{3}) \times
\dot{H}^{s_{p}-1}(\mathbb{R}^{3})} + \| D^{s_{p}
-\frac{1}{2}}(|v_{[1]}|^{p-1}v_{[1]})  \|_{ \tilde{W} (I_{j+1}) } \\
& \lesssim \| ( v_{[1]}(t_{j}),\partial_{t} v_{[1]}(t_{j}) )
\|_{\dot{H}^{s_{p}}(\mathbb{R}^{3}) \times
\dot{H}^{s_{p}-1}(\mathbb{R}^{3})}  + \| D^{s_{p} -\frac{1}{2}}
v_{[1]}  \|_{W(I_{j+1})} \| v_{[1]} \|^{p-1}_{S(I_{j+1})} \\
& \lesssim \| v_{[1]} \|_{X(I_{j})} + \eta^{p-1} \| v_{[1]}
\|_{X(I_{j+1})}
\end{array}
\end{equation}
Notice that $l \lesssim 1$: this follows from Conjecture
\ref{Conj:1}, Condition \ref{Conj:2} and the inequality

\begin{equation}
\begin{array}{ll}
\| (u_{0},u_{1}) \|_{\dot{H}^{s_{p}}(\mathbb{R}^{3}) \times
\dot{H}^{s_{p}-1}(\mathbb{R}^{3})} & \leq \sup_{t \in I_{max,g_{1}}}
\|(u(t), \partial_{t}u(t) ) \|_{\dot{H}^{s_{p}}(\mathbb{R}^{3})
\times \dot{H}^{s_{p}-1}(\mathbb{R}^{3})} \\  & \leq
C_{2}\left( \| (u_{0},u_{1}) \|_{\tilde{H}^{2} \times \tilde{H}^{1}} \right) \\
& \lesssim 1
\end{array}
\label{Eqn:ControlHsp}
\end{equation}
,following from Condition \ref{Conj:2} and the assumption \\ $\|(
u_{0}, u_{1} ) \|_{\tilde{H}^{2} \times \tilde{H}^{1}} \leq 1$
\footnote{Notice that, at this stage, we only need to know that $ \|
( u_{0}, u_{1}) \|_{\dot{H}^{s_{p}}(\mathbb{R}^{3}) \times
\dot{H}^{s_{p}-1}(\mathbb{R}^{3}) } \leq \| ( u_{0}, u_{1})
\|_{\tilde{H}^{2} \times \tilde{H}^{1} } \leq 1$ and apply
Conjecture \ref{Conj:1}. Therefore the introduction of $ \sup_{t \in
I_{max,g_{1}}} \| (u(t),\partial_{t} u(t))
\|_{\dot{H}^{s_{p}}(\mathbb{R}^{3}) \times
\dot{H}^{s_{p}-1}(\mathbb{R}^{3})}$ in (\ref{Eqn:ControlHsp}) in
redundant. This is done on purpose. Indeed, we will use Condition
\ref{Conj:2} in other parts of the argument: see
(\ref{Eqn:ControlHsp1}) }. Therefore, by a standard continuity
argument and iteration on $j$ we have

\begin{equation}
\begin{array}{ll}
\| v_{[1]}  \|_{X(\mathbb{R})} & \lesssim 1
\end{array}
\label{Eqn:Controlv1Y}
\end{equation}

\item \textit{Substep 2}: we control $ \|u_{[1]} - v_{[1]} \|_{X([-\tilde{t},\tilde{t}])}
$ for $\tilde{t} <<  1$ to be chosen later. By time reversal
symmetry it is enough to control $\| u_{[1]} - v_{[1]}
\|_{X([0,\tilde{t}])}$. To this end we consider $w_{[1]} := u_{[1]}
- v_{[1]}$. We get

\begin{equation}
\begin{array}{ll}
\partial_{tt} w_{[1]} - \triangle w_{[1]} = & -|w_{[1]} +
v_{[1]}|^{p-1} (v_{[1]} + w_{[1]}) g_{1}(v_{[1]} + w_{[1]}) + 2
|v_{[1]}|^{p-1} v_{[1]}
\end{array}
\end{equation}
Let $\eta^{'}<< 1$. By (\ref{Eqn:Controlv1Y}), we can divide $[0,
\tilde{t}]$ into subintervals $(J_{k})_{1 \leq k \leq m }$ that
satisfy the following properties

\begin{enumerate}

\item $ \| D^{s_{p}-\frac{1}{2}} v_{[1]} \|_{L_{t}^{\infty -}
L_{x}^{3+}(J_{k})} = \eta^{'} $ or $ \| D^{s_{p} -\frac{1}{2}}
v_{[1]} \|_{W(J_{k})} = \eta^{'} $ for $ 1 \leq k < m$

\item $ \| D^{s_{p} - \frac{1}{2}} v_{[1]} \|_{W(J_{k})} \leq
\eta^{'} $ and $ \|  D^{s_{p} - \frac{1}{2}} v_{[1]}
\|_{L_{t}^{\infty -} L_{x}^{3+} (J_{k})  } \leq \eta^{'}$ for $1
\leq k \leq m$

\end{enumerate}
We have

\begin{equation}
\begin{array}{ll}
\| w_{[1]} \|_{X(J_{k+1})} & \lesssim \| ( w_{[1]}(t_{k}),
\partial_{t} w_{[1]}(t_{k}) \|_{\dot{H}^{s_{p}}(\mathbb{R}^{3})
\times \dot{H}^{s_{p}-1}(\mathbb{R}^{3})}  \\
& + \| D^{s_{p}-\frac{1}{2}} ( 2 |v_{[1]}|^{p-1} v_{[1]} - 2
|v_{[1]} +
w_{[1]}|^{p-1} (v_{[1]} +w_{[1]}) )  \|_{\tilde{W}(J_{k+1})}  \\
& + \| D^{s_{p}-\frac{1}{2}} ( h_{1}(|v_{[1]}+w_{[1]}|)
|v_{[1]}+w_{[1]}|^{p-1} (v_{[1]}+w_{[1]}) \|_{L_{t}^{1}
L_{x}^{\frac{3}{2}}(J_{k+1})}
\end{array}
\end{equation}
Let

\begin{equation}
\begin{array}{ll}
A_{1} & := \| D^{s_{p}-\frac{1}{2}} ( 2 |v_{[1]}|^{p-1} v_{[1]} - 2
|v_{[1]} + w_{[1]}|^{p-1} (v_{[1]} +w_{[1]}) )
\|_{\tilde{W}(J_{k+1})}
\\
A_{2} & : = \| D^{s_{p}-\frac{1}{2}} ( h_{1}(|v_{[1]}+w_{[1]}|)
|v_{[1]}+w_{[1]}|^{p-1} (v_{[1]}+w_{[1]}) \|_{L_{t}^{1}
L_{x}^{\frac{3}{2}}(J_{k+1})}
\end{array}
\end{equation}
By  the fractional Leibnitz rule applied to $q(x):=|x|^{p-1} x
h(x)$, (\ref{Eqn:hpne}), (\ref{Eqn:hprimeonecd}), Sobolev embedding
and H\"older in time we have

\begin{equation}
\begin{array}{ll}
A_{2} & \lesssim \| |v_{[1]} + w_{[1]}|^{\frac{p-1}{2}+}
\|_{L_{t}^{1+} L_{x}^{3-}
(J_{k+1})} \| D^{s_{p}-\frac{1}{2}} (v_{[1]}+w_{[1]}) \|_{L_{t}^{\infty -} L_{x}^{3 +}(J_{k+1})} \\
& \lesssim \| v_{[1]} + w_{[1]}
\|^{\frac{p-1}{2}+}_{L_{t}^{\frac{p-1}{2}+}
L_{x}^{\frac{3(p-1)+}{2}} (J_{k+1})} \| D^{s_{p}-\frac{1}{2}}
(v_{[1]}+w_{[1]}) \|_{L_{t}^{\infty -} L_{x}^{3 +}(J_{k+1})} \\
& \lesssim \tilde{t}  \| D^{s_{p} -\frac{1}{2}} (v_{[1]} +
w_{[1]}) \|^{\frac{p+1}{2} +}_{L_{t}^{\infty -} L_{x}^{3 +}(J_{k+1})} \\
& \lesssim \tilde{t}   \eta^{\frac{p+1}{2}+} + \tilde{t} \|
D^{s_{p}-\frac{1}{2}} w_{[1]}
\|^{\frac{p+1}{2}+}_{L_{t}^{\infty -} L_{x}^{3 +}(J_{k+1}) } \\
& \lesssim \tilde{t}   \eta^{\frac{p+1}{2}+} + \tilde{t} \| w_{[1]}
\|^{\frac{p+1}{2}+}_{X(J_{k+1})}
\end{array}
\end{equation}

As for $A_{1}$ we follow \cite{kenigmerle1},p 9

\begin{equation}
\begin{array}{ll}
A_{1} & \lesssim \left( \| v_{[1]} \|^{p-1}_{S(J_{k+1})} + \|
w_{[1]} \|^{p-1}_{S(J_{k+1})} \right) \| D^{s_{p} -\frac{1}{2}}
w_{[1]} \|_{W(J_{k+1})} + \left( \| v_{[1]} \|^{p-2}_{S(J_{k+1})} +
\| w_{[1]}
\|^{p-2}_{S(J_{k+1})} \right) \\
&   \left( \| D^{s_{p} - \frac{1}{2}} v_{[1]} \|_{W(J_{k+1})} + \|
D^{s_{p} - \frac{1}{2}}w_{[1]} \|_{W(J_{k+1})} \right) \|
w_{[1]} \|_{S(J_{k+1})} \\ \\
& \lesssim \eta^{p-1} \| w_{[1]} \|_{X(J_{k+1})} + \| w_{[1]}
\|^{p}_{X(J_{k+1})} + \eta^{p-2} \| w_{[1]} \|_{X(J_{k+1})}^{2} +
\eta \| w_{[1]} \|_{X(J_{k+1})}^{p-1}
\end{array}
\label{Eqn:XJkInduc}
\end{equation}
This follows from (\ref{Eqn:FracLeibnDiff}) and
(\ref{Eqn:SobEmbed2}). Therefore we have

\begin{equation}
\begin{array}{ll}
\| w_{[1]} \|_{X(J_{k+1})} & \lesssim  \| w_{[1]} \|_{X(J_{k})} +
\eta^{\frac{p+1}{2}+} \tilde{t} + \tilde{t} \| w_{[1]}
\|^{\frac{p+1}{2}+}_{X(J_{k+1})} + \eta^{p-1} \|
w_{[1]} \|_{X(J_{k+1})} + \| w_{[1]} \|^{p}_{X(J_{k+1})} \\
& + \eta^{p-2} \| w_{[1]} \|^{2}_{Y(J_{k+1})} + \eta \| w_{[1]}
\|^{p-1}_{X(J_{k+1})}
\end{array}
\label{Eqn:Eq1}
\end{equation}
Let $C$ be the constant determined by (\ref{Eqn:Eq1}). By induction,
we have

\begin{equation}
\begin{array}{ll}
\| w_{[1]} \|_{X(J_{k})} & \leq (2C)^{k} \tilde{t}
\end{array}
\end{equation}
provided that for $1 \leq k \leq m-1$ we have

\begin{equation}
\begin{array}{ll}
C \eta^{\frac{p+1}{2}} \tilde{t} & << C (2C)^{k}
\tilde{t} \\
C \tilde{t} \left( (2C)^{k} \tilde{t} \right)^{\frac{p+1}{2}+} &
<< (2C)^{k} \tilde{t} \\
C \eta^{p-1} (2C)^{k+1} \tilde{t} & << C (2C)^{k}
\tilde{t} \\
C \left( (2C)^{k} \tilde{t} \right)^{p} & << C (2C)^{k}
\tilde{t} \\
C \eta^{p-2} \left( (2C)^{k+1} \tilde{t} \right)^{2}
& << C (2C)^{k} \tilde{t} \\
\eta \left( (2C)^{k+1} \tilde{t} \right)^{p-1} & << C (2C)^{k}
\tilde{t}
\end{array}
\end{equation}
These inequalities are satisfied if $\eta << 1$ and

\begin{equation}
\begin{array}{ll}
\tilde{t} & <<  1
\end{array}
\end{equation}
since $k \leq m-1$ and, by (\ref{Eqn:Controlv1Y}), $m \lesssim 1$.
We conclude that

\begin{equation}
\begin{array}{ll}
\| w_{[1]} \|_{X([0,\tilde{t}])} & \lesssim  1 \\
\end{array}
\end{equation}

\item \textit{Substep 3}: we control $\| u_{[1]} \|_{X([-T,T])}$. By time reversal symmetry, it is
enough to control $\| u_{[1]} \|_{X([0,T])}$. Recall that $T \leq
1$. We chop $T \leq 1$ into subintervals $(J_{k^{'}}=[a_{k^{'}},
b_{k^{'}}] )_{1 \leq k^{'} \leq l^{'}}$ such that $| J_{k^{'}} | =
\tilde{t}$ for $1 \leq k^{'} < l^{'}$ and $ |J_{l^{'}}| \leq
\tilde{t}$. Notice that, by Condition \ref{Conj:2}, we have

\begin{equation}
\begin{array}{ll}
\| ( u(a_{k^{'}}), \partial_{t} u (a_{k^{'}}))
\|_{\dot{H}^{s_{p}}(\mathbb{R}^{3}) \times
\dot{H}^{s_{p}-1}(\mathbb{R}^{3})} & \leq \sup_{t \in I_{max,g_{1}}}
\|(u(t), \partial_{t}u(t) ) \|_{\dot{H}^{s_{p}}(\mathbb{R}^{3})
\times \dot{H}^{s_{p}-1}(\mathbb{R}^{3})} \\  & \leq
C_{2} \left( \| (u_{0}, u_{1}) \|_{\tilde{H}^{2} \times \tilde{H}^{1}} \right) \\
& \lesssim 1
\end{array}
\label{Eqn:ControlHsp1}
\end{equation}
,taking advantage of the assumption $\|( u_{0}, u_{1} )
\|_{\tilde{H}^{2} \times \tilde{H}^{1}} \leq 1$. For each $k^{'}$ we
define $v_{[1,k^{'}]}$ to be the solution of \footnote{ in
particular $v_{[1,k^{'}]} = v_{[1]}$}

\begin{equation}
\left\{
\begin{array}{ll}
\partial_{tt} v_{[1,k^{'}]} -\triangle v_{[1,k^{'}]} & = -
|v_{[1,k^{'}]}|^{p-1} v_{[1,k^{'}]}   \\
 v_{[1,k^{'}]}(a_{k^{'}}) & = u_{[1]}(a_{k^{'}}) \\
\partial_{t} v_{[1,k^{'}]}(a_{k^{'}})  & = \partial_{t} u_{[1]}(a_{k^{'}})
\end{array}
\right.
\end{equation}
By slightly modifying the proof of the previous step and letting
$v_{[1,k^{'}]}$ play the role of $v_{[1]}$ (see previous step), this
leads, by (\ref{Eqn:ControlHsp1}), to

\begin{equation}
\begin{array}{ll}
\| v_{[1,k^{'}]} \|_{X(\mathbb{R})} & \lesssim 1
\end{array}
\end{equation}
and

\begin{equation}
\begin{array}{ll}
\| w_{[1,k^{'}]} \|_{X(J_{k^{'}})} & \lesssim  1
\end{array}
\end{equation}
with $ w_{[1,k^{'}]}: =u_{[1]} - v_{[1,k^{'}]}$. Therefore $ \|
u_{[1]} \|_{ X(J_{k^{'}})} \lesssim 1$ and summing over the
$J_{k^{'}}$ s we have

\begin{equation}
\begin{array}{ll}
\| u_{[1]} \|_{X([0,T])} & \lesssim  1
\end{array}
\label{Eqn:Controlw1critical}
\end{equation}

\textit{Substep $4$}: we control $ \| ( u_{[1]}, \partial_{t}
u_{[1]} ) \|_{ L_{t}^{\infty} \tilde{H}^{2} ([-1,1]) \times
L_{t}^{\infty} \tilde{H}^{1} ([-1,1])}$ and  $ \| u_{[1]}
\|_{S([-1,1])} $. We get from (\ref{Eqn:Controlw1critical})

\begin{equation}
\begin{array}{ll}
\| u_{[1]} \|_{S([-1,1])} & \lesssim  1
\end{array}
\label{Eqn:Controlu1S}
\end{equation}
\end{itemize}
By Proposition \ref{Prop:controlLtinfty} and (\ref{Eqn:Controlu1S})
we have

\begin{equation}
\begin{array}{ll}
\| ( u_{[1]}, \partial_{t} u_{[1]} ) \|_{ L_{t}^{\infty}
\tilde{H}^{2} ([-1,1]) \times L_{t}^{\infty} \tilde{H}^{1} ([-1,1])}
& \lesssim 1
\end{array}
\end{equation}
Therefore

\begin{equation}
\begin{array}{ll}
\max{ \left( \| u_{[1]} \|_{S([-1,1])}, \, \| ( u_{[1]},
\partial_{t} u_{[1]} ) \|_{L_{t}^{\infty} \tilde{H}^{2} ([-1,1]) \times
L_{t}^{\infty} \tilde{H}^{1} ([-1,1])}   \right)}  & \lesssim 1
\end{array}
\label{Eqn:Controlw1critical2}
\end{equation}
We let $C^{'}_{1}$ (defined in the statement of Lemma \ref{lem:gi})
be equal to one. We can assume without the loss of generality that
the constant determined by $\lesssim$ in
(\ref{Eqn:Controlw1critical2}) is larger than $1$. We let $C_{1}$
(defined in the statement of Lemma \ref{lem:gi}) be equal to this
constant. $C^{'}_{1}$ and $C_{1}$ satisfy (\ref{Eqn:Ci1}) and
(\ref{Eqn:Ci2}).

\textit{Step $2$}: construction of $g_{i}$ from $g_{i-1}$

Recall that $[-T,T] \subset [-i,i]$ and  $\| (u_{0},u_{1})
\|_{\tilde{H}^{2} \times \tilde{H}^{1}} \leq i $.  From
(\ref{Eqn:Constructiongi}) it is enough to construct $g_{i}$ for
$|x| > A C_{i-1}$. It is clear that, by choosing $C^{'}_{i}$ large
enough, we can construct find a function $\tilde{g}_{i}$ defined on
$[A C_{i-1}, \, C^{'}_{i}]$ such that $g_{i}$, defined by

\begin{equation}
\begin{array}{ll}
g_{i}(x) & : = \left\{
\begin{array}{l}
g_{i-1}(x), \, |x| \leq A C_{i-1}  \\
\tilde{g}_{i}(x), \, C^{'}_{i}  \geq |x| \geq A C_{i-1} \\
i+1, \, |x| \geq C^{'}_{i}
\end{array}
\right.
\end{array}
\end{equation}
is smooth, slowly increasing; it satisfies (\ref{Eqn:Conditiong2}),
(\ref{Eqn:Conditiong3}), and

\begin{equation}
\begin{array}{ll}
\int_{A C_{i-1}} ^{C^{'}_{i}} \frac{1}{y g_{i}^{2}(y) } \, dy & \geq
i
\end{array}
\end{equation}
It remains to determine $C_{i}$ (defined in the statement of Lemma
\ref{lem:gi}). To do that we slightly modify the step $1$.

We sketch the argument. Let $h_{i}(x):= g_{i}(x) - (i+1)$. Then
$h_{i}(x)=0$ if $|x|
> C^{'}_{i}$. It is not difficult to see that

\begin{equation}
\begin{array}{ll}
|h_{i}(x)| & \lesssim_{i} \frac{1}{|x|^{\frac{p-1}{2}+}}
\end{array}
\label{Eqn:Boundhi}
\end{equation}

\begin{equation}
\begin{array}{ll}
|h^{'}_{i}(x)| & \lesssim_{i} \frac{1}{|x|^{\frac{p+1}{2}+}}
\end{array}
\label{Eqn:Boundhprimei}
\end{equation}
Let $u_{[i]}$ and $v_{[i]}$ be the solutions of the equations

\begin{equation}
\left\{
\begin{array}{ll}
\partial_{tt} u_{[i]} - \triangle u_{[i]} & = - |u_{[i]}|^{p-1}
u_{[i]} g_{i}(|u_{[i]}|) \\
u_{[i]}(0) & := u_{0}  \\
\partial_{t} u_{[i]}(0) & := u_{1}
\end{array}
\right.
\end{equation}
and

\begin{equation}
\left\{
\begin{array}{ll}
\partial_{tt} v_{[i]} - \triangle v_{[i]} & = - (i+1) |v_{[i]}|^{p-1}
v_{[i]} \\
v_{[i]}(0) & := u_{0} \\
\partial_{t} v_{[i]}(0) & := u_{1}
\end{array}
\right. \label{Eqn:visol}
\end{equation}
We have

\begin{itemize}

\item \textit{Substep 1}: we have

\begin{equation}
\begin{array}{ll}
\| v_{[i]} \|_{X(\mathbb{R})} & \lesssim_{i} 1
\end{array}
\end{equation}
, by adapting the proof of Substep $1$, Step $2$. Notice, in
particular, that we can use Conjecture $1$ and control $ \| v_{[i]}
\|_{S(\mathbb{R})} $ since $w_{[i]}:= (i+1)^{\frac{1}{p-1}} v_{[i]}
$ satisfies $ \partial_{tt} w_{[i]} - \triangle w_{[i]} = -
|w_{[i]}|^{p-1} w_{[i]}$.

\item \textit{Substep 2}: we have  $ \| u_{[i]} - v_{[i]} \|_{X([0,\tilde{t}])} \lesssim_{i} 1$
for $\tilde{t} <<_{i} 1$, by adapting the proof of Substep $2$, Step
$2$. The dependance on $i$ basically comes from (\ref{Eqn:Boundhi}),
(\ref{Eqn:Boundhprimei}) and (\ref{Eqn:visol}).

\item \textit{Substep 3}: we prove that $\| u_{[i]} \|_{X([-T,T])} \lesssim_{i,p}
1$. By time reversal symmetry, it is enough control $\| u_{[i]}
\|_{X([0,T])}$. Recall that $T \leq i$. We chop $[0,T]$ into
subintervals $(J_{k^{'}}=[a_{k^{'}}, b_{k^{'}}] )_{1 \leq k^{'} \leq
l^{'}}$ such that $| J_{k^{'}} | = \tilde{t}$ for $1 \leq k^{'} <
l^{'}$ and $ |J_{l^{'}}| \leq \tilde{t}$ (with $\tilde{t}$ defined
in Substep $2$). By Condition \ref{Conj:2} and the assumption $\| (
u_{0},u_{1}) \|_{\tilde{H}^{2} \times \tilde{H}^{1}} \leq i$, we
have

\begin{equation}
\begin{array}{ll}
\| ( u_{[i]}(a_{k^{'}}), \partial_{t} u_{[i]} (a_{k^{'}}))
\|_{\dot{H}^{s_{p}}(\mathbb{R}^{3}) \times
\dot{H}^{s_{p}-1}(\mathbb{R}^{3})} & \leq \sup_{t \in
I_{\max,g_{i}}} \left\| ( u_{[i]}(t), \, \partial_{t} u_{[i]}(t) )
\right\|_{\dot{H}^{s_{p}}(\mathbb{R}^{3})
\times \dot{H}^{s_{p}-1}(\mathbb{R}^{3})} \\
& \lesssim_{i} 1
\end{array}
\label{Eqn:Controlnormi}
\end{equation}
We introduce

\begin{equation}
\left\{
\begin{array}{ll}
\partial_{tt} v_{[i,k^{'}]} -\triangle v_{[i,k^{'}]} & = -(i+1)
|v_{[i,k^{'}]}|^{p-1} v_{[i,k^{'}]}   \\
 v_{[i,k^{'}]}(a_{k^{'}}) & = u_{[i]}(a_{k^{'}}) \\
\partial_{t} v_{[i,k^{'}]}(a_{k^{'}})  & = \partial_{t} u_{[i]}(a_{k^{'}})
\end{array}
\right.
\end{equation}
and, by using (\ref{Eqn:Controlnormi}), we can prove

\begin{equation}
\begin{array}{ll}
\| u_{[i]} \|_{S([-i,i])} & \lesssim_{i} 1
\end{array}
\label{Eqn:ControlCritNormi}
\end{equation}

\textit{Substep $4$}. By using Proposition \ref{Prop:controlLtinfty}
and (\ref{Eqn:ControlCritNormi}) we get

\begin{equation}
\begin{array}{ll}
\max{ \left( \| u_{[i]} \|_{S([-i,i])}, \, \| ( u_{[i]},
\partial_{t} u_{[i]} ) \|_{L_{t}^{\infty} \tilde{H}^{2} ([-i,i]) \times
L_{t}^{\infty} \tilde{H}^{1} ([-i,i])}   \right)}  & \lesssim_{i} 1
\end{array}
\label{Eqn:Controlw1criticali}
\end{equation}
We can assume without loss of generality that the constant
determined by $\lesssim$ in (\ref{Eqn:Controlw1criticali}) is larger
than $i$ and $C^{'}_{i}$. We let $C_{i}$ be equal to this constant.
(\ref{Eqn:Ci1}) and (\ref{Eqn:Ci2}) are satisfied.

\end{itemize}
This ends the proof.
\end{itemize}

\end{document}